\documentclass[a4paper,american,reqno]{amsart}

\usepackage{babel}
\usepackage[utf8]{inputenc}
\usepackage[T1]{fontenc}
\usepackage{lmodern}
\usepackage{amssymb}
\usepackage{graphicx}
\usepackage{xcolor}
\usepackage{booktabs}
\usepackage{array}
\usepackage{xspace}
\usepackage{microtype}
\usepackage{etoolbox}
\usepackage{csquotes}
\usepackage{tikz}
\usetikzlibrary{arrows.meta,positioning,calc}
\usepackage{caption}
\usepackage{subcaption}
\usepackage{placeins}

\usepackage[style=authoryear,
  dashed=false,
  giveninits=true,
  uniquename=init,
  maxbibnames=100,
  maxcitenames=2,
  isbn=false,
  natbib=true,
  backend=bibtex]{biblatex}
\usepackage[colorlinks,citecolor=blue,urlcolor=blue,linkcolor=blue]{hyperref}
\usepackage{cleveref}

\usepackage[hmargin=27mm,vmargin=25mm]{geometry}
\graphicspath{{./}}
\DeclareCaptionFont{wiley}{\normalsize}
\AtBeginEnvironment{tabular}{\normalsize}

\newtheorem{theorem}{Theorem}

\newtheorem{corollary}[theorem]{Corollary}
\newtheorem{proposition}[theorem]{Proposition}
\newtheorem{assumption}{Assumption}
\theoremstyle{definition}
\newtheorem{definition}[theorem]{Definition}
\newtheorem{example}[theorem]{Example}
\newtheorem{remark}[theorem]{Remark}

\usepackage{pgfplots}
\pgfplotsset{compat=1.18}

\newcommand{\define}{\mathrel{{\mathop:}{=}}}
\newcommand{\set}[1]{\{#1\}}
\newcommand{\Set}[1]{\left\{#1\right\}}
\newcommand{\abs}[1]{\lvert #1\rvert}
\newcommand{\card}{\abs}
\newcommand{\Inedges}{\delta^{-}}
\newcommand{\Outedges}{\delta^{+}}

\newcommand{\Z}{\mathbb{Z}}

\newcommand{\Defset}[3][\defsep]{\Set{#2#1#3}}

\newenvironment{widetable}[1][htb]
  {\begin{table}[#1]}{\end{table}}
\newcommand{\PenaltyPanelWidth}{0.31\textwidth}
\newcommand{\PenaltyPanelBreak}{\hfill}

\definecolor{exblue}{HTML}{0072B2}
\definecolor{exorange}{HTML}{E69F00}
\definecolor{exgreen}{HTML}{009E73}
\definecolor{exred}{HTML}{D55E00}
\tikzset{
  vtx/.style   = {circle, draw, thick, fill=white, minimum size=7.5mm,
                  inner sep=0pt, font=\normalsize},
  arc/.style   = {->, >={Stealth[length=2.4mm]}, semithick,
                  shorten >=1pt, black!55},
  len/.style   = {fill=white, inner sep=1.2pt, font=\small, text=black},
  tree/.style        = {draw=exblue, very thick},
  longest/.style     = {draw=exorange, very thick},
  detour/.style      = {draw=exgreen, very thick},
  flowone/.style     = {draw=exblue, very thick},
  flowtwo/.style     = {draw=exorange, very thick},
  flowshared/.style  = {draw=exgreen, very thick},
  interdicted/.style = {draw=exred, very thick, dashed},
  root/.style        = {draw=exblue, very thick, circle, minimum size=9mm},
  ex/12/.style={}, ex/21/.style={}, ex/15/.style={}, ex/25/.style={},
  ex/32/.style={}, ex/34/.style={}, ex/43/.style={}, ex/41/.style={},
  ex/53/.style={}, ex/54/.style={},
  showcosts/.code = {\excoststrue},
}
\newif\ifexcosts
\newcommand{\ExLabel}[2]{\ifexcosts{\scriptsize$(#1,#2)$}\else$#1$\fi}
\newcommand{\FlowLabel}[2]{\scriptsize$(#1,#2)$}

\newcommand{\RunningExample}{%
  \node[vtx] (v1) at (0,3.2)   {$1$};
  \node[vtx] (v2) at (3.2,3.2) {$2$};
  \node[vtx] (v3) at (3.2,0)   {$3$};
  \node[vtx] (v4) at (0,0)     {$4$};
  \node[vtx] (v5) at (1.6,1.6) {$5$};
  \draw[arc,ex/12] (v1) to[bend left=18] node[len,above] {\ExLabel{3}{1}} (v2);
  \draw[arc,ex/21] (v2) to[bend left=18] node[len,below] {\ExLabel{3}{2}} (v1);
  \draw[arc,ex/43] (v4) to[bend left=18] node[len,above] {\ExLabel{5}{2}} (v3);
  \draw[arc,ex/34] (v3) to[bend left=18] node[len,below] {\ExLabel{1}{1}} (v4);
  \draw[arc,ex/32] (v3) -- node[len,right] {\ExLabel{2}{2}} (v2);
  \draw[arc,ex/41] (v4) -- node[len,left] {\ExLabel{3}{1}} (v1);
  \draw[arc,ex/15] (v1) -- node[len,pos=0.6] {\ExLabel{3}{2}} (v5);
  \draw[arc,ex/25] (v2) -- node[len,pos=0.62] {\ExLabel{1}{1}} (v5);
  \draw[arc,ex/53] (v5) -- node[len,pos=0.55] {\ExLabel{4}{1}} (v3);
  \draw[arc,ex/54] (v5) -- node[len,pos=0.55] {\ExLabel{3}{2}} (v4);
}

\newcommand{\SurvivingFlowGraph}{%
  \node[vtx] (v1) at (0,3.2)   {$1$};
  \node[vtx] (v2) at (3.2,3.2) {$2$};
  \node[vtx] (v3) at (3.2,0)   {$3$};
  \node[vtx] (v4) at (0,0)     {$4$};
  \node[vtx] (v5) at (1.6,1.6) {$5$};
  \draw[arc] (v2) to[bend left=18] node[len,below] {\FlowLabel{3}{2}} (v1);
  \draw[arc,flowtwo] (v4) to[bend left=18] node[len,above] {\FlowLabel{5}{2}} (v3);
  \draw[arc] (v3) to[bend left=18] node[len,below] {\FlowLabel{1}{1}} (v4);
  \draw[arc,flowshared] (v3) -- node[len,right] {\FlowLabel{2}{2}} (v2);
  \draw[arc] (v4) -- node[len,left] {\FlowLabel{3}{1}} (v1);
  \draw[arc,flowshared] (v1) -- node[len,pos=0.6] {\FlowLabel{3}{2}} (v5);
  \draw[arc] (v2) -- node[len,pos=0.62] {\FlowLabel{1}{1}} (v5);
  \draw[arc,flowone] (v5) -- node[len,pos=0.55] {\FlowLabel{4}{1}} (v3);
  \draw[arc,flowtwo] (v5) -- node[len,pos=0.55] {\FlowLabel{3}{2}} (v4);
}

\newif\ifpreprintwidetable
\AtBeginEnvironment{widetable}{\preprintwidetabletrue}
\newsavebox{\preprinttablebox}
\BeforeBeginEnvironment{tabular}{%
  \ifpreprintwidetable\begin{lrbox}{\preprinttablebox}\fi}
\AfterEndEnvironment{tabular}{%
  \ifpreprintwidetable
    \end{lrbox}%
    \ifdim\wd\preprinttablebox>\linewidth
      \resizebox{\linewidth}{!}{\usebox{\preprinttablebox}}%
    \else
      \usebox{\preprinttablebox}%
    \fi
  \fi}
\tikzset{x=1.2cm,y=1.2cm,len/.append style={font=\normalsize},
  every node/.append style={transform shape=false}}
\renewcommand{\ExLabel}[2]{\ifexcosts$(#1,#2)$\else$#1$\fi}
\renewcommand{\FlowLabel}[2]{$(#1,#2)$}

\let\preprintabstract\abstract
\let\preprintendabstract\endabstract
\RenewDocumentCommand{\abstract}{O{ABSTRACT}+m}{%
  \preprintabstract#2\preprintendabstract}
\let\preprintkeywords\keywords
\newcommand{\preprintsubjectclass}[1]{\subjclass[2020]{#1}}
\ExplSyntaxOn
\RenewDocumentCommand{\keywords}{m}{
  \tl_set:Nn \l_tmpa_tl {#1}
  \tl_replace_all:Nnn \l_tmpa_tl {~|~} {,~}
  \exp_args:NV \preprintkeywords \l_tmpa_tl
}
\NewDocumentCommand{\MSC}{m}{
  \tl_set:Nn \l_tmpa_tl {#1}
  \tl_replace_all:Nnn \l_tmpa_tl {;~} {,~}
  \exp_args:NV \preprintsubjectclass \l_tmpa_tl
}
\ExplSyntaxOff
\let\bmsubsection\section

\begin{document}

\title[Multiple-Pairs Shortest-Path Interdiction]{An Effective MILP
  Formulation for the\\Multiple-Pairs Shortest-Path Interdiction
  Problem\\via Strengthened Penalties}

\author[A. L. Croella, F. Furini, S. Martin, M. Schmidt, J. Th\"urauf]{Anna Livia
  Croella, Fabio Furini, Sébastien Martin,\\Martin Schmidt, Johannes Th\"urauf}

\address[A. L. Croella]{Department of Engineering and Science, Universitas Mercatorum,
Piazza Mattei 10, 00186, Rome, Italy}
\email{annalivia.croella@unimercatorum.it}

\address[F. Furini]{Department of Computer, Control and Management
  Engineering ``Antonio Ruberti'', Sapienza University of Rome,
  Via Ariosto 25, 00185 Roma, Italy}
\email{fabio.furini@uniroma1.it}

\address[S. Martin]{Laboratoire de Conception, Optimisation et Mod\'elisation des Syst\`emes (LCOMS), Universit\'e de Lorraine, France}
\email{sebastien.martin@univ-lorraine.fr}

\address[M. Schmidt]{Department of Mathematics, Trier University,
  Universit\"atsring 15, 54296 Trier, Germany}
\email{martin.schmidt@uni-trier.de}

\address[J. Th\"urauf]{Department of Liberal Arts and Social Sciences,
  Discrete Optimization, University of Technology Nuremberg,
  Nuremberg, Germany}
\email{johannes.thuerauf@utn.de}

\date{\today}

\keywords{Network interdiction | Shortest path problems | Multiple-pairs
  shortest paths | Bilevel optimization | Big-$M$ values | Mixed-integer
  linear optimization | Combinatorial optimization}

\MSC{90C10; 90C35; 90C57}

\abstract[ABSTRACT]{%
  We study the interdiction of the multiple-pairs shortest-path problem:
  a leader removes arcs of a directed graph, within a budget on their
  interdiction costs, to maximize the sum of the shortest-path distances
  of a prescribed set of origin-destination pairs, all of which remain
  connected after every feasible removal. Single-pair and all-pairs
  shortest-path interdiction are special cases. After modeling arc
  removal by penalties added to the lengths of the interdicted arcs,
  linear optimization duality yields a compact mixed-integer
  linear formulation whose strength depends on these penalties. We characterize exactly
  which penalties are valid, and we prove that, for every set of pairs,
  the valid penalties admit a componentwise smallest vector, which gives
  the strongest formulation within this family of coefficients. Although
  the interdicted arcs are removed simultaneously, each entry of this
  vector is determined by residual shortest-path distances along paths
  that contain the penalized arc as their only interdicted arc. We then
  derive penalties computable in polynomial time, which exploit the
  budget, the interdiction costs, and the requested destinations through
  families of paths of which at least one survives every feasible
  interdiction. {With a time limit of $1800$~seconds,
  the formulation using the strongest of the proposed
  poly-time computable penalty rules solves $190$ of the $229$ all-pairs
  instances, compared with $52$ using a uniform penalty. The linear
  optimization relaxation is solved for every rule on $196$ instances.
  On this common subset, the strongest rule reduces the mean relaxation
  gap from $102.5\%$ to $4.6\%$.}
  Experiments with sets of pairs show that the destination-aware
  penalties strengthen the relaxation further.
}

\maketitle

\section{Introduction}
\label{intro}

Let $\mathcal G=(\mathcal V,\mathcal A)$ be a directed graph in which
every arc~$a\in\mathcal A$ has a positive integer length~$\ell_a$, and
let
$\mathcal K\subseteq\set{(s,t)\in\mathcal V\times\mathcal V:s\neq t}$
be a prescribed set of ordered origin-destination pairs, called the
\emph{requested pairs}. Given $(\mathcal G,{\boldsymbol\ell},\mathcal K)$,
the \emph{multiple-pairs shortest-path problem} (MPSP) consists of
determining, for every requested pair~$(s,t)\in\mathcal K$, a shortest
$s$--$t$ path and its length~$d(s,t)$, called the \emph{distance}
from~$s$ to~$t$.
We call the sum $\sum_{(s,t)\in\mathcal K}d(s,t)$ of the distances of
the requested pairs the \emph{MPSP value}. The \emph{single-pair
shortest-path problem} is the special case $\card{\mathcal K}=1$ and the
\emph{all-pairs shortest-path problem} (APSP) is the special case in
which $\mathcal K$ contains every ordered pair of distinct vertices.

In this paper, we study the interdiction counterpart of the MPSP,
which we call the \emph{MPSP interdiction problem}. Every arc~$a\in\mathcal A$ also has a positive integer interdiction cost~$c_a$,
and a \emph{leader} removes a set~$I\subseteq\mathcal A$ of arcs whose
total cost does not exceed a budget~$b\in\Z_{\ge0}$. The remaining
graph~$\mathcal G(\mathcal V,\mathcal A\setminus I)$ is called the
\emph{interdicted network}. The \emph{follower} then routes every requested
pair along a shortest path of the interdicted network, and the leader
maximizes the sum of the resulting distances, i.e., it solves
\[
  \max\Defset{
    \sum_{(s,t)\in\mathcal K}d^I(s,t)
  }{
    I\subseteq\mathcal A,\ \textstyle\sum_{a\in I}c_a\le b
  },
\]
where $d^I(s,t)$ denotes the distance from~$s$ to~$t$ in the interdicted
  network, called the residual distance. Thus, the leader maximizes the MPSP value of the
interdicted network. We assume that every feasible removal leaves a directed
path for every requested pair, so that this sum is finite. The
classic \emph{shortest-path interdiction problem}, with a single pair,
and \emph{APSP interdiction}, with all pairs, are special cases. The reduction
of~\citet[Section~5]{Khachiyan_etal:2008} shows that shortest-path
interdiction is NP-hard also under this assumption
(Remark~\ref{rem:hardness} in \Cref{sec:problem-statement}). Hence, MPSP
interdiction is NP-hard already with a single requested
pair.

Network interdiction has been studied for many follower problems; see
\citet{Wood:1993} and \citet{Smith_Song:2020} for overviews. In APSP
interdiction, the objective is the sum of all pairwise distances of the
interdicted network, i.e., the directed counterpart of the Wiener
index~\citep{Wiener:1947}. This objective links the problem to
\emph{distance-based critical node and edge detection}, in which
vertices or edges of a (usually undirected) graph are removed to degrade
the distances between its vertices
\citep{Veremyev_etal:2015,Hooshmand_etal:2020,Alozie_etal:2024}, and
to the single most vital edge problem with respect to all pairs of
\citet{Venema_etal:2000}, in which one edge is removed.

A standard way to model arc removal in shortest-path interdiction is to
add a penalty to the length of every interdicted arc instead of
deleting it~\citep{Israeli_Wood:2002}. The feasible region of the
follower then no longer depends on the decision of the leader, and
linear optimization duality turns the bilevel problem into a
single-level mixed-integer linear optimization problem, in which the binary
interdiction variables are multiplied by penalties (big-$M$ values). A
penalty that is too small lets the follower use an interdicted arc and
underestimates the distances, while a penalty that is too large weakens
the linear optimization relaxation. Indeed, for a requested
origin~$w$, a fractional value~$x_a$ in the relaxation adds $x_aM(a,w)$ to the length of arc~$a$: with
a large penalty, a small fraction of the budget spread over many arcs
already lengthens all of them substantially, so that the relaxation
overestimates the residual distances and its bound can be far from the
optimal value. The quality of the penalties therefore determines the
strength of the formulation, and with it the size of the
branch-and-bound trees and the instances that can be solved to proven
optimality. For arc destruction,
\citet{Israeli_Wood:2002} use a uniform penalty, an artificial delay
added to every interdicted arc, and observe that both the direct
formulation and their Benders decomposition, whose cuts also contain
the penalty, become weak and slow as the penalty grows, whereas a
penalty that is too small yields invalid solutions. They therefore
recommend a covering decomposition whose master problem contains no
penalty terms. Its covering constraints are supervalid inequalities
without lifting. More generally, choosing correct
big-$M$ values is known to be a delicate issue in bilevel
optimization~\citep{Kleinert_etal:2020,Pineda_Morales:2019}.

Structured penalties are not new in network interdiction. For
multicommodity flow interdiction, \citet{Lim_Smith:2007} derive arc-
and commodity-dependent penalties computable in polynomial time, and
distinguish penalties that make the follower avoid the interdicted arcs
from penalties that only preserve its optimal value. Their validity
conditions are sufficient, and their follower may leave demand unmet,
so that the penalty of an arc only has to cancel the reward collected
along the paths through that arc. In MPSP interdiction, instead, every
requested pair must remain connected. The penalties are valid only if
every path through interdicted arcs is at least as long as the
residual shortest path after the \emph{joint} removal of all
interdicted arcs. The required penalties are therefore determined
jointly by residual detours, which depend on the other interdicted arcs
and, through them, on the interdiction costs and the budget. This dependence is specific to
distance objectives, and it is the subject of this paper.

\subsection{Contribution and Outline}
\label{sec:contribution}

The main contribution of this paper is a theory of valid penalties for
MPSP interdiction, together with penalties computable in polynomial
time. To the best of our knowledge, neither an exact characterization
of valid penalties nor penalties that exploit the interdiction costs and
the budget have been studied for shortest-path interdiction. Our
contributions are threefold:

{
\begin{enumerate}
\item[(i)]
We characterize the valid penalties exactly and prove that,
for every pair set, they admit a componentwise smallest vector.
This vector yields the strongest linear optimization relaxation among
all valid choices of penalties for our compact formulation.
Within a single interdiction, the penalties of different arcs can
compensate each other. Over the family of feasible interdictions, which
is closed under taking subsets, this coupling disappears, and each
smallest penalty is determined by paths that contain the penalized arc
as their only interdicted arc and end at a requested destination. In
APSP interdiction, these values can be expressed directly in terms of the
residual distances between the tail and the head of each arc. We
compute them
by enumeration on small instances and use them only as a theoretical reference.

\item[(ii)] We derive penalties computable in polynomial time from upper
bounds on residual distances. The arc-length bounds sum the largest arc
lengths, and the surviving-path bounds come from minimum-cost flows
whose paths cannot all be destroyed by a feasible interdiction, so that
they exploit the interdiction costs and the budget. The root-dependent
versions subtract the distance from the root to the tail of the arc,
which every path through the arc already pays. The destination-aware
versions, in addition, impose the conditions only for the requested
destinations and use the distance from the head of the arc to each of
them. Together, the two refinements are the counterpart, for
distances, of the prefix--suffix bounds of \citet{Lim_Smith:2007}.

\item[(iii)] We assess computationally the strength of the linear optimization
relaxation, the cost of computing the penalties, and the performance of
the compact formulation, against a uniform penalty of the type used
by~\citet{Israeli_Wood:2002}, and against the covering decomposition and
the supervalid inequalities of the same authors. Separately, on small
all-pairs instances, we compare the penalties of the poly-time computable rules
entry by entry with the smallest valid ones. With the strongest rule,
the compact formulation solves $190$ of $229$ all-pairs instances within
$1800$ seconds, against $52$ with the uniform penalty, and this rule
attains the smallest penalty on $86.6\%$ of the entries of the small
instances, and on $76.7\%$ of the entries on which the smallest
penalty is positive. On the $124$ all-pairs instances on which the covering
decomposition is tested, the compact formulation solves $116$, against
$13$ for the covering decomposition, which remains the fastest method on
the single-pair grids of~\citet{Israeli_Wood:2002}. On pair sets, the
destination-aware penalties strengthen the relaxation when few
destinations are requested, and the covering decomposition solves fewer
instances as the number of pairs grows.
\end{enumerate}
}

{
The remainder of the paper is organized as follows.
\Cref{sec:literature} reviews the related literature.
\Cref{sec:problem-statement} defines MPSP interdiction and its natural
and penalty formulations. \Cref{sec:approaches} presents the compact
formulation, which extends the dualization of~\citet{Israeli_Wood:2002}
to several pairs and serves as a framework. \Cref{sec:bigM-analysis}
contains the theory of the penalties and the poly-time computable rules, and
\Cref{sec:svi} adapts the supervalid inequalities
of~\citet{Israeli_Wood:2002} to the compact formulation.
\Cref{sec:computational} reports the computational experiments, and
\Cref{sec:conclusions} concludes the paper.}

\section{Related Literature}
\label{sec:literature}

In this section, we review three streams of literature: shortest-path
interdiction and most vital arcs, distance-based critical node and edge
detection, and the choice of penalties (big-$M$ values) in interdiction
and bilevel reformulations. For general overviews of network
interdiction and of mixed-integer bilevel optimization, we refer
to~\citet{Smith_Song:2020} and~\citet{Kleinert_etal:2021}, respectively.

\subsection{Shortest-Path Interdiction and Most Vital Arcs}

\citet{Fulkerson_Harding:1977} and~\citet{Golden:1978} study the
interdiction of a single shortest $s$--$t$ path when arcs can be
lengthened by continuous amounts, a setting solvable by minimum-cost
flow techniques. With binary interdiction decisions, the problem becomes
NP-hard~\citep{Khachiyan_etal:2008}. \citet{Israeli_Wood:2002} study the
discrete problem on directed graphs, in which an interdicted arc~$a$ is
lengthened by a delay~$d_a$. They dualize the shortest-path follower,
which gives a compact single-level formulation (the
dualize-and-combine approach in the terminology
of~\citet{Smith_Song:2020}), and propose a Benders decomposition
strengthened by supervalid inequalities. Arc
destruction is modeled by a sufficiently large artificial delay~$d^+$,
e.g., $d^+=\abs{N}\max_{a}c_a$ in their notation. Israeli and
Wood point out that both the compact formulation and the Benders
decomposition suffer from weak LP relaxations when the delays are large,
whereas a delay that is too small yields invalid solutions. For arc destruction they
therefore recommend a covering decomposition whose master problem
contains no delay terms. The cardinality version, i.e., the $k$
\emph{most vital arcs} problem, goes back to~\citet{Malik_etal:1989}. Inapproximability results are given by~\citet{Khachiyan_etal:2008}, and
a parameterized complexity analysis for undirected graphs
by~\citet{Bazgan_etal:2019}. Further variants include randomized
interdiction strategies~\citep{Holzmann_Smith:2021}, fortification of
arcs before the attack, solved by \citet{Lozano_Smith:2017} by sampling
$s$--$t$ paths of the follower, and the blocker
version, in which a minimum number of vertices is deleted so that no
short $s$--$t$ path survives~\citep{Magnouche20}. All these works
consider a single origin-destination pair. With respect to all pairs,
\citet{Venema_etal:2000} determine the single edge of a weighted
undirected graph whose removal increases the sum of all pairwise
distances the most.

\citet{Ayyildiz_etal:2018} study single-source, multiple-sink
shortest-path interdiction on undirected networks using an aggregated
flow, a uniform interdiction penalty, and LP dualization, anticipating
the one-root formulation used here. \citet{Ayyildiz_etal:2023} consider
a capacitated multiple-source, multiple-sink extension in which the
follower chooses which sources supply the sinks. Source capacities
couple these choices, unlike the independent routing of our prescribed
pairs. Neither study provides the exact penalty characterization or
the bounds from surviving paths developed here.

Our follower evaluates every pair in~$\mathcal K$ under the same
interdiction, and we model arc destruction through penalties, as in the
compact formulation of~\citet{Israeli_Wood:2002}. Instead of avoiding the
penalties by a decomposition, we study which penalty values keep the
compact formulation exact and how small they can be chosen.

\subsection{Distance-Based Critical Node and Edge Detection}

Critical node and edge detection problems delete vertices or edges,
within a budget, to degrade a connectivity measure of the interdicted network.
\citet{Veremyev_etal:2014} propose an integer programming framework with
joint node and edge deletions, for measures based on the sizes of the
connected components. \citet{Veremyev_etal:2015} introduce the
distance-based critical node detection problem in undirected graphs,
whose objective is a weighted sum, over the vertex pairs, of a function
of their residual distance. One of its classes maximizes the Wiener index,
with a large constant assigned to disconnected pairs, and the edge-weighted
case is also considered. \citet{Aringhieri_etal:2019} study the
complexity of this problem on paths, trees, and series-parallel graphs. For the Wiener-index class, it is strongly NP-hard on trees with
arbitrary pair weights, and they give a pseudo-polynomial dynamic program
on trees with unit pair weights and unit deletion costs.
For this class, \citet{Hooshmand_etal:2020}
propose a formulation with pair-distance variables, bounded through the
edge lengths and big-$M$ terms and linked by triangle inequalities, and
solve it by a Benders decomposition in which the dual variables of the
subproblems are computed analytically. Their constant~$M$ also gives the contribution of disconnected pairs.
Path-based integer programming formulations for other classes of this
problem are proposed by~\citet{Alozie_etal:2021a}, and exact methods for
minimizing the number of vertex pairs connected by a path of length at
most~$k$ by~\citet{Salemi_Buchanan:2022}.
Edge deletions are considered by~\citet{Veremyev_etal:2019}, who degrade
the closeness centrality of a given group of vertices in undirected
unweighted graphs with edge removal costs, and
by~\citet{Alozie_etal:2024}, who reduce the edge version to the node
version on a graph with subdivided edges and also give a direct edge
formulation, in which distances above a threshold are replaced by a
constant.

The problem studied here removes arcs of a directed graph, with arbitrary
interdiction costs, and its objective is the sum of the untruncated
distances of the ordered pairs in~$\mathcal K$, all of which remain
connected under every feasible interdiction. Our formulation has one
block of distance labels per active root, in which the penalties lengthen
the interdicted arcs. The pair-distance and triangle-inequality
formulation of~\citet{Hooshmand_etal:2020} is designed for a different
problem, the removal of a bounded number of vertices of an undirected
graph, in which disconnected pairs are allowed and contribute a
constant to the objective. Consider its directed adaptation to arc
interdiction without the disconnection caps, with a distance
variable~$\theta_{ij}$ for each ordered pair~$(i,j)$, a binary
interdiction variable~$x_a$ for each arc, the arc constraints
$\theta_{uv}\le\ell_a+M(a)x_a$ for $a=(u,v)$, and the triangle
inequalities, where $M(a)$ is a root-independent penalty. For every
fixed value of the interdiction variables, also fractional, the
shortest-path distances with lengths~$\ell_a+M(a)x_a$ satisfy all these
constraints, and applying the triangle inequalities along a shortest
path bounds every distance variable by the corresponding distance. The distance labels of the compact formulation of this paper take the same
values for the requested pairs. Hence, with the same root-independent penalties,
the adaptation has the same linear optimization bound as our
formulation. Root-dependent penalties cannot be inserted in its arc
constraints, which are shared by all roots, and the adaptation is not
used as a benchmark in this paper.

\subsection{Penalties in Interdiction and Bilevel Reformulations}

The closest predecessor of our penalty analysis is the multicommodity
flow interdiction study of~\citet{Lim_Smith:2007}. In their penalty
formulation, the interdiction of arc~$a$ reduces the reward of
commodity~$k$ on that arc by a penalty~$M_{ak}$, and the follower's
problem is dualized to obtain a single-level formulation. They derive three
levels of penalties: a uniform value, one value per commodity, and one
value per arc and commodity, obtained from longest paths that can be
computed in polynomial time because the reward graphs have no positive
cycles. Penalties strictly larger than these bounds force the follower to
avoid interdicted arcs, while penalties equal to them only preserve the
optimal value. Both conditions are sufficient, and they rely on dummy
arcs that allow unmet demand, so that the flow on an interdicted arc can
always be dropped. In their experiments on discrete interdiction, the
arc-and-commodity penalties reduce the running time of the penalty
formulation only slightly on one test set, but by more than one order of
magnitude on the hardest instances of another. In hub interdiction, \citet{Ramamoorthy_etal:2018} dualize a min-cost
routing follower whose paths visit at most two hubs, linearize the
products of the dual variables and the hub variables with big-$M$
constraints, and derive a valid value of each~$M_k$ from the difference
between the costs of the most and the least expensive routing path of
every origin-destination pair. Since the routing paths are polynomially
many, they also replace the follower by closest-assignment constraints,
which require no dual variables.

In the survey of~\citet[Section~6.2]{Kleinert_etal:2021}, the penalty
approach to discrete interdiction requires penalties large enough that
the follower does not use interdicted objects, and the choice of these
coefficients is described as crucial for the efficiency of the
single-level reformulation.
In bilevel optimization, big-$M$ values are mostly discussed for the
reformulation of the lower-level optimality conditions.
\citet{Kleinert_etal:2020} prove that verifying that a big-$M$ bounding
the lower-level dual variables cuts off no optimal dual solution is as
hard as solving the bilevel problem itself, and~\citet{Pineda_Morales:2019}
show by a counterexample that big-$M$ values tuned by trial and error can
lead to highly suboptimal solutions. Our penalties lengthen arcs in the
follower's objective and do not bound multipliers, so these hardness
results do not apply directly. They motivate, however, an exact notion of
validity. For optimal transmission switching, \citet{Fattahi_etal:2019}
show that the strongest big-$M$ values are NP-hard to compute and derive
bounds from shortest paths in a fixed connected spanning subnetwork of
lines that cannot be switched off.

With respect to~\citet{Lim_Smith:2007}, the dualized penalty formulation,
penalties depending on the arc and on the commodity (here, the active
root), and the distinction between strict and value-preserving validity
are not new. What changes with distance objectives is that a requested
pair cannot be dropped: the penalties must be large enough that no path
through interdicted arcs is shorter than the residual shortest path
after their joint removal, and therefore depend on the budget and on the
interdiction costs. In contrast to the fixed spanning subnetwork
of~\citet{Fattahi_etal:2019}, only the arcs with $c_a>b$ are protected
from interdiction in our setting, and they need not contain a connected
spanning subnetwork. The surviving paths depend on the interdiction. Our contribution lies in the necessary and sufficient
conditions for valid penalties, which yield the componentwise smallest
penalties for every pair set, the all-pairs case included, and in poly-time computable strengthened
penalties built from families of paths of which at least one survives
every feasible interdiction. Our destination-aware rules are the distance
analogue of the prefix and suffix bounds of~\citet{Lim_Smith:2007}.

\section{Problem Statement and Formulations}
\label{sec:problem-statement}

Let $\mathcal{G}=(\mathcal{V},\mathcal{A})$ be a strongly connected
directed graph with $n$ vertices, $\mathcal{V}=\set{1,\dots,n}$, and
$m\define\card{\mathcal{A}}$ arcs. The graph has no loops and no
parallel arcs. We write $a=(u,v)$ for an arc, with tail~$u$ and head~$v$.
Each arc has positive integer length~$\ell_a$. For
$I\subseteq\mathcal A$, $\mathcal G(\mathcal V,\mathcal A\setminus I)$
is the interdicted network, i.e., the graph obtained from~$\mathcal G$ by
removing the arcs of~$I$. We denote by $d^I(u,v)$ the length of a shortest
$u$--$v$ path in the interdicted network, with $d^I(u,v)=+\infty$ if no such
path exists, and set~$d(u,v)\define d^\emptyset(u,v)$.

For the pair set~$\mathcal K \subseteq \Set{(w,t)\in\mathcal{V} \times \mathcal{V}: w\neq t}$, define the active roots
\[
  \mathcal R\define
  \Set{w\in\mathcal V:\text{ there exists }t\text{ with }(w,t)\in\mathcal K},
\]
i.e., the vertices that are the origin of at least one requested pair,
and, for every $w\in\mathcal R$, the requested destinations
\[
  \mathcal T_w\define\Set{t\in\mathcal V:(w,t)\in\mathcal K},
\]
i.e., the vertices that form a requested pair with origin~$w$.
The MPSP value is
\[
  \varphi_{\mathcal K}(\mathcal G,{\boldsymbol \ell})
  \define
  \sum_{w\in\mathcal R}\sum_{t\in\mathcal T_w}d(w,t),
\]
i.e., the sum of the shortest-path distances of the requested pairs.
The all-pairs specialization is
\[
  \mathcal K^{\mathrm{AP}}
  \define
  \Set{(w,t)\in\mathcal V\times\mathcal V:w\neq t},
\]
i.e., the set of all ordered pairs of distinct vertices.
Then, $\mathcal R=\mathcal V$ and
$\mathcal T_w=\mathcal V\setminus\set{w}$. In this specialization we
write $\varphi(\mathcal G,{\boldsymbol \ell})$ instead of
$\varphi_{\mathcal K^{\mathrm{AP}}}(\mathcal G,{\boldsymbol \ell})$. The single-pair specialization is
$\mathcal K=\set{(s,t)}$, with $\mathcal R=\set{s}$ and
$\mathcal T_s=\set{t}$. Its interdiction counterpart, defined below,
is the shortest-path interdiction problem of~\citet{Israeli_Wood:2002}
with arc removal.

\begin{figure}[htb]
  \centering
  \begin{subfigure}[t]{0.44\textwidth}\centering\vspace{0pt}
    \begin{tikzpicture}[ex/21/.style=tree, ex/25/.style=tree,
                        ex/53/.style=tree, ex/54/.style=tree]
      \RunningExample
      \node[root] at (v2) {};
    \end{tikzpicture}
    \caption{Graph, arc lengths~$\ell_a$, and the shortest paths from
      vertex~$2$.}
    \label{FIG_a_graph}
  \end{subfigure}\hfill
  \begin{subfigure}[t]{0.52\textwidth}\centering\vspace{0pt}
    \renewcommand{\arraystretch}{1.25}
    \begin{tabular}{@{}c|ccccc|c@{}}
      from\,$\backslash$\,to & $1$ & $2$ & $3$ & $4$ & $5$ & sum\\
      \hline
      $1$ & -- & $3$ & $7$ & $6$ & $3$ & $19$\\
      \textcolor{exblue}{$2$} & \textcolor{exblue}{$3$} & \textcolor{exblue}{--}
        & \textcolor{exblue}{$5$} & \textcolor{exblue}{$4$}
        & \textcolor{exblue}{$1$} & \textcolor{exblue}{$13$}\\
      $3$ & $4$ & $2$ & -- & $1$ & $3$ & $10$\\
      $4$ & $3$ & $6$ & $5$ & -- & $6$ & $20$\\
      $5$ & $6$ & $6$ & $4$ & $3$ & -- & $19$\\
      \hline
      \multicolumn{6}{@{}r|}{$\varphi(\mathcal{G},{\boldsymbol \ell})$}
        & $\mathbf{81}$
    \end{tabular}
    \medskip
    \caption{Shortest-path distances.}
    \label{FIG_a_matrix}
  \end{subfigure}
  \caption{The running example, shown in the all-pairs specialization: a strongly connected digraph with
    $n=5$ vertices and $m=10$ arcs (left), and the lengths of the
    shortest paths between all ordered pairs of vertices (right). In
    blue, the shortest paths from vertex~$2$ (circled) and the
    corresponding row of distances.}
  \label{FIG_a}
\end{figure}

All figures and numerical examples of the running example, except
Example~\ref{ex:running-pairs}, use
$\mathcal K=\mathcal K^{\mathrm{AP}}$, while the formulations and
theoretical results below are stated for a general pair set whenever
possible.

Each arc also has a positive integer interdiction cost~$c_a$ and the
leader has budget~$b$. Define
\[
  \psi_{\mathcal K}(\mathcal G,{\boldsymbol \ell},{\boldsymbol c},b)
  \define
  \max\Defset{
    \varphi_{\mathcal K}(\mathcal G(\mathcal V,\mathcal A\setminus I),
                         {\boldsymbol \ell})
  }{
    I\subseteq\mathcal A,\
    \textstyle\sum_{a\in I}c_a\le b
  },
\]
i.e., the largest MPSP value of an interdicted network obtained by
removing a set of arcs whose total interdiction cost is at most the
budget.
For $\mathcal K=\mathcal K^{\mathrm{AP}}$ we keep the shorter notation
$\psi(\mathcal G,{\boldsymbol \ell},{\boldsymbol c},b)$.

\subsection{The Multiple-Pairs Shortest-Path Problem}
\label{sec:all-pairs-shortest}

As already mentioned above, we make the following first
  assumption.

\begin{assumption}
  \label{ass:connected}
  The graph~$\mathcal G=(\mathcal V,\mathcal A)$ is directed and
    strongly connected.
\end{assumption}

For $u\in\mathcal V$, the sets of outgoing and incoming arcs are denoted
by $\Outedges(u)$ and $\Inedges(u)$. Let $\mathcal K$ be the prescribed
pair set, $\mathcal R$ the active roots, and $\mathcal T_w$ the requested
destinations of root~$w$.

For every $w\in\mathcal R$ and $a\in\mathcal A$, let $y_{aw}$ denote
the flow of root~$w$ on arc~$a$. An integer formulation of the MPSP is
\begin{subequations}
  \label{APSP-PRIMAL}
  \begin{alignat}{3}
    \varphi_{\mathcal K}(\mathcal G,{\boldsymbol \ell}) =\ 
       &\makebox[2em][c]{$\displaystyle\min_{\boldsymbol y}$}
       &\quad&\sum_{w\in\mathcal R}\sum_{a\in\mathcal A}\ell_a y_{aw}
       &&
       \label{OBJ_APSP_PRIMAL}\\*[1.2ex]
       &\makebox[2em][c]{s.t.}
       &&\sum_{a\in\Inedges(t)}y_{aw}
       -\sum_{a\in\Outedges(t)}y_{aw}
       \ge1,
       &\quad&w\in\mathcal R,\ t\in\mathcal T_w,
       \label{C1_APSP_PRIMAL}\\[1.2ex]
    &&&\sum_{a\in\Inedges(v)}y_{aw}
       -\sum_{a\in\Outedges(v)}y_{aw}
       \ge0,
       &\quad&w\in\mathcal R,\
         v\in\mathcal V\setminus(\mathcal T_w\cup\Set{w}),
       \label{C0_MPSP_PRIMAL}\\[1.2ex]
    &&&y_{aw}\in\mathbb Z_{\ge0},
       &\quad&w\in\mathcal R,\ a\in\mathcal A.
       \label{C2_APSP_PRIMAL}
  \end{alignat}
\end{subequations}

The objective function~\eqref{OBJ_APSP_PRIMAL} minimizes the total
length of the flows of all active roots. The \emph{demand
constraints}~\eqref{C1_APSP_PRIMAL} require every requested
destination~$t\in\mathcal T_w$ to receive at least one unit of the flow
of root~$w$, and the \emph{transit constraints}~\eqref{C0_MPSP_PRIMAL}
prevent the flow from being created at the other non-root vertices. Constraints~\eqref{C2_APSP_PRIMAL} define the domain of the variables.
Hence, for every active root $w$, one unit is absorbed by each requested
destination and every other non-root vertex is a transit vertex. Since
all lengths are positive, an optimal solution contains no circulation
and decomposes into $\card{\mathcal T_w}$ simple requested paths.

\begin{remark}[Total unimodularity]
  \label{rem:TU}
  For every root~$w$, the constraint matrix is the vertex-arc
  incidence matrix of~$\mathcal G$ without the row of~$w$. Hence, the
  complete matrix is block diagonal and totally unimodular.
\end{remark}

We can therefore replace integrality by
\[
  y_{aw}\ge0,\quad w\in\mathcal R,\ a\in\mathcal A.
\]
Using a nonnegative dual variable~$\pi_{vw}$
(also called a node potential) for every
$w\in\mathcal R$ and $v\neq w$, with $\pi_{ww}\define0$, gives
\begin{subequations}
  \label{AP-PATH-DUAL}
  \begin{alignat}{3}
    \varphi_{\mathcal K}(\mathcal G,{\boldsymbol \ell}) =\ 
       &\makebox[2em][c]{$\displaystyle\max_{\boldsymbol \pi}$}
       &\quad&\sum_{w\in\mathcal R}\sum_{t\in\mathcal T_w}\pi_{tw}
       &&
       \label{OBJ_D_AP}\\*[1.2ex]
       &\makebox[2em][c]{s.t.}
       &&\pi_{vw}-\pi_{uw}\le\ell_{(u,v)},
       &\quad&w\in\mathcal R,\
        (u,v)\in\mathcal A\setminus\Outedges(w),
       \label{C1_D_AP}\\[1.2ex]
    &&&\pi_{vw}\le\ell_{(w,v)},
       &\quad&w\in\mathcal R,\
        (w,v)\in\Outedges(w),
       \label{C2_D_AP}\\[1.2ex]
    &&&\pi_{vw}\ge0,
       &\quad&w\in\mathcal R,\ v\in\mathcal V\setminus\Set{w}.
       \label{C3_D_AP}
  \end{alignat}
\end{subequations}

The variable~$\pi_{vw}$ can be read as a distance label from root~$w$
to vertex~$v$. The objective function~\eqref{OBJ_D_AP} maximizes the
sum of the labels of the requested destinations. The
constraints~\eqref{C1_D_AP} require that the label of the head of every
arc exceed the label of its tail by at most the length of the arc, and
Constraints~\eqref{C2_D_AP} impose the same condition on the arcs
leaving the root, whose label is zero. Constraints~\eqref{C3_D_AP}
define the domain of the variables. Hence, every label is at most the
length of every path from the root to its vertex.
For every requested destination~$t\in\mathcal T_w$, an optimal dual
solution satisfies $\pi_{tw}=d(w,t)$. For non-destination vertices, the
optimal label need not equal the shortest-path distance because its
objective coefficient is zero. This additional dual freedom is the key
difference from the all-pairs case.

When $\mathcal K=\mathcal K^{\mathrm{AP}}$,
$\mathcal T_w=\mathcal V\setminus\set{w}$ holds and the zero-demand family
disappears. Hence, the two formulations become the APSP formulations.

\subsection{MPSP Interdiction}

\subsubsection{A Natural Formulation}

Let $x_a=1$ if and only if arc~$a$ is interdicted. The leader's feasible set is
\[
  \mathcal X\define
    \Set{{\boldsymbol x}\in\Set{0,1}^{\mathcal A}:
      \sum_{a\in\mathcal A}c_ax_a\le b}.
\]

For every active
root~$w$, let $\kappa_w\define\card{\mathcal T_w}$, which bounds the
flow of root~$w$ on every arc in an optimal follower solution.
By Remark~\ref{rem:TU} applied to the interdicted network,
whenever the requested pairs are connected, the follower has an
integral optimal solution. Since all lengths are positive, this
solution decomposes into $\kappa_w$ simple paths for each root~$w$,
one to each destination in~$\mathcal T_w$, and uses every arc at most
$\kappa_w$ times. Hence, bounding each arc flow of root~$w$
by~$\kappa_w$ cuts off no optimal follower solution.
The direct bilevel formulation is
\begin{subequations}
  \label{eq:apspip}
  \begin{alignat}{3}
    \psi_{\mathcal K}(\mathcal G,{\boldsymbol \ell},{\boldsymbol c},b)
       = \max_{{\boldsymbol x}\in\mathcal X}\quad
       &\makebox[2em][c]{$\displaystyle\min_{\boldsymbol y}$}
       &\quad&\sum_{w\in\mathcal R}\sum_{a\in\mathcal A}\ell_a y_{aw}
       &&
       \label{OBJ_apspip}\\*[1.2ex]
       &\makebox[2em][c]{s.t.}
       &&\sum_{a\in\Inedges(t)}y_{aw}
       -\sum_{a\in\Outedges(t)}y_{aw}
       \ge1,
       &\quad&w\in\mathcal R,\ t\in\mathcal T_w,
       \label{C_flow_apspip}\\[1.2ex]
    &&&\sum_{a\in\Inedges(v)}y_{aw}
       -\sum_{a\in\Outedges(v)}y_{aw}
       \ge0,
       &\quad&w\in\mathcal R,\
        v\in\mathcal V\setminus(\mathcal T_w\cup\Set{w}),
       \label{C_zero_apspip}\\[1.2ex]
    &&&y_{aw}\le\kappa_w(1-x_a),
       &\quad&w\in\mathcal R,\ a\in\mathcal A,
       \label{C_bigM_apspip}\\[1.2ex]
    &&&y_{aw}\ge0,
       &\quad&w\in\mathcal R,\ a\in\mathcal A.
       \label{C_dom_y_apspip}
  \end{alignat}
\end{subequations}

The leader chooses the interdiction decision~${\boldsymbol x}$ to
maximize the objective function~\eqref{OBJ_apspip}, which is the
optimal value of the follower, i.e., the total length of the shortest
paths of the requested pairs in the interdicted network. The follower
solves the inner minimization problem of
Formulation~\eqref{eq:apspip} for the given~${\boldsymbol x}$:
Constraints~\eqref{C_flow_apspip} and~\eqref{C_zero_apspip} are its
demand and transit constraints, and the \emph{linking
constraints}~\eqref{C_bigM_apspip} forbid any flow on an interdicted
arc. The set~$\mathcal X$ imposes the leader's budget and binary
restrictions. Constraints~\eqref{C_dom_y_apspip} define the domains of
the follower variables.

Finite MPSP objective values only require the requested pairs to remain
connected.

\begin{assumption}
  \label{ass:feasibility}
  For every feasible interdiction decision and every
  $(w,t)\in\mathcal K$, the interdicted network contains a directed
  $w$--$t$ path.
\end{assumption}

The following proposition gives an equivalent condition for
Assumption~\ref{ass:feasibility} that can be verified using
minimum-cut computations.
\begin{proposition}
  \label{obs:cut}
  Assumption~\ref{ass:feasibility} holds if and only if, for every
  $(w,t)\in\mathcal K$, every $(w,t)$-cut has total interdiction cost
  larger than~$b$.
\end{proposition}
\begin{proof}
  A feasible interdiction disconnects a requested pair if and only if
  it contains a complete cut separating that pair.
\end{proof}

Thus, the condition can be checked with at most
$\card{\mathcal K}$ minimum-cut computations.

Some results below concern distances between vertices that need not
form a requested pair, for instance all the pairs in the all-pairs
case. For those results we use the following stronger assumption.

\begin{assumption}
  \label{ass:strong-residual}
  Every feasible interdiction leaves the interdicted network strongly
  connected.
\end{assumption}

Assumption~\ref{ass:strong-residual} implies
Assumption~\ref{ass:feasibility} for every pair set~$\mathcal K$, and
the two coincide for $\mathcal K=\mathcal K^{\mathrm{AP}}$.

\begin{remark}
  \label{rem:single-arc}
  Since interdiction costs are positive integers, a feasible
  interdiction set contains at most~$b$ arcs. For fixed~$b$, the MPSP
  interdiction problem is solved in polynomial time by enumerating the
  $O(m^b)$ feasible interdiction sets and evaluating
  $\varphi_{\mathcal K}$ on each interdicted network.
\end{remark}

\begin{remark}
  \label{rem:hardness}
  Shortest-path interdiction remains NP-hard when every
  feasible interdiction leaves the origin and the destination connected,
  and even when it leaves the graph strongly connected.
  \citet[Section~5]{Khachiyan_etal:2008} reduce the minimum
  vertex cover problem on a graph with $n$ vertices to the decision
  version of the problem (can at most~$k$ removable arcs raise the
  $s$--$t$ distance to at least $\lambda=n-1$?) on a directed
  graph~$\mathcal H$ with arc lengths in~$\set{0,1}$, whose
  non-removable arcs contain an $s$--$t$ path~$P_0$; we may assume
  $1\le k\le n$. In our setting, the non-removable arcs receive
  cost~$k+1$, the removable arcs cost~$1$, and $b=k$, so that the
  feasible interdictions are exactly the sets of at most~$k$ removable
  arcs. Let $N$ be the number of vertices of~$\mathcal H$. Replacing
  every length~$\ell_a$ by $N\ell_a+1$ and the threshold~$\lambda$
  by~$N\lambda$ preserves the answer: the new length of every path is at
  least $N$ times its old length plus one, and some shortest path for
  the old lengths is simple, with fewer than~$N$ arcs, so that an
  integer distance~$d$ becomes a value in $[Nd+1,Nd+N-1]$. Let
  $L$ be the new length of~$P_0$, so that $d^I(s,t)\le L$ for every
  feasible~$I$. Add a vertex~$z$ and, for every vertex~$v$
  of~$\mathcal H$, the arcs $(v,z)$ of length~$L$ and $(z,v)$ of
  length~$1$, both of cost~$k+1$. The graph remains simple, every
  interdicted network is strongly connected, since every vertex reaches
  every other one through~$z$, and every $s$--$t$ path through~$z$ has
  length at least $L+1>d^I(s,t)$, so that no residual $s$--$t$ distance
  changes. All the numbers of the construction are polynomial in~$N$,
  so that the problem is NP-hard in the strong sense.
\end{remark}

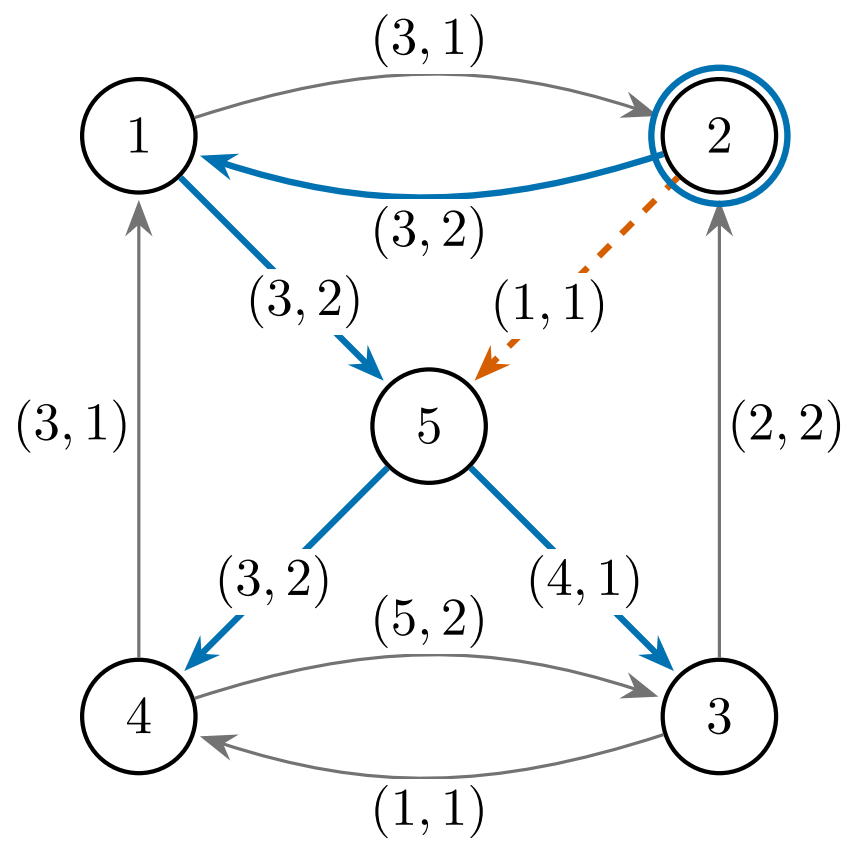
\begin{figure}[htb]
  \centering
  \newcommand{\chg}[1]{\textcolor{exred}{\boldmath$\mathbf{#1}$}}
  \renewcommand{\arraystretch}{1.2}
  \begin{subfigure}[t]{0.44\textwidth}\centering\vspace{0pt}
    \begin{tikzpicture}[showcosts, ex/25/.style=interdicted, ex/21/.style=tree,
                        ex/15/.style=tree, ex/53/.style=tree,
                        ex/54/.style=tree]
      \RunningExample
      \node[root] at (v2) {};
    \end{tikzpicture}
    \caption{$b=1$: arc~$(2,5)$ is interdicted, and the new shortest paths from vertex~$2$ are shown.}
    \label{FIG_solution_b1}
  \end{subfigure}\hfill
  \begin{subfigure}[t]{0.52\textwidth}\centering\vspace{0pt}
    \begin{tabular}{@{}c|ccccc|c@{}}
      from\,$\backslash$\,to & $1$ & $2$ & $3$ & $4$ & $5$ & sum\\
      \hline
      $1$ & -- & $3$ & $7$ & $6$ & $3$ & $19$\\
      $2$ & $3$ & -- & \chg{10} & \chg{9} & \chg{6} & \chg{28}\\
      $3$ & $4$ & $2$ & -- & $1$ & \chg{7} & \chg{14}\\
      $4$ & $3$ & $6$ & $5$ & -- & $6$ & $20$\\
      $5$ & $6$ & $6$ & $4$ & $3$ & -- & $19$\\
      \hline
      \multicolumn{6}{@{}r|}{$\psi(\mathcal{G},{\boldsymbol \ell},{\boldsymbol c},1)$} & \chg{100}
    \end{tabular}
    \medskip
    \caption{$b=1$: distances in $\mathcal{G}(\mathcal{V},\mathcal{A}\setminus I)$, $I=\set{(2,5)}$.}
    \label{FIG_solution_b1_matrix}
  \end{subfigure}

  
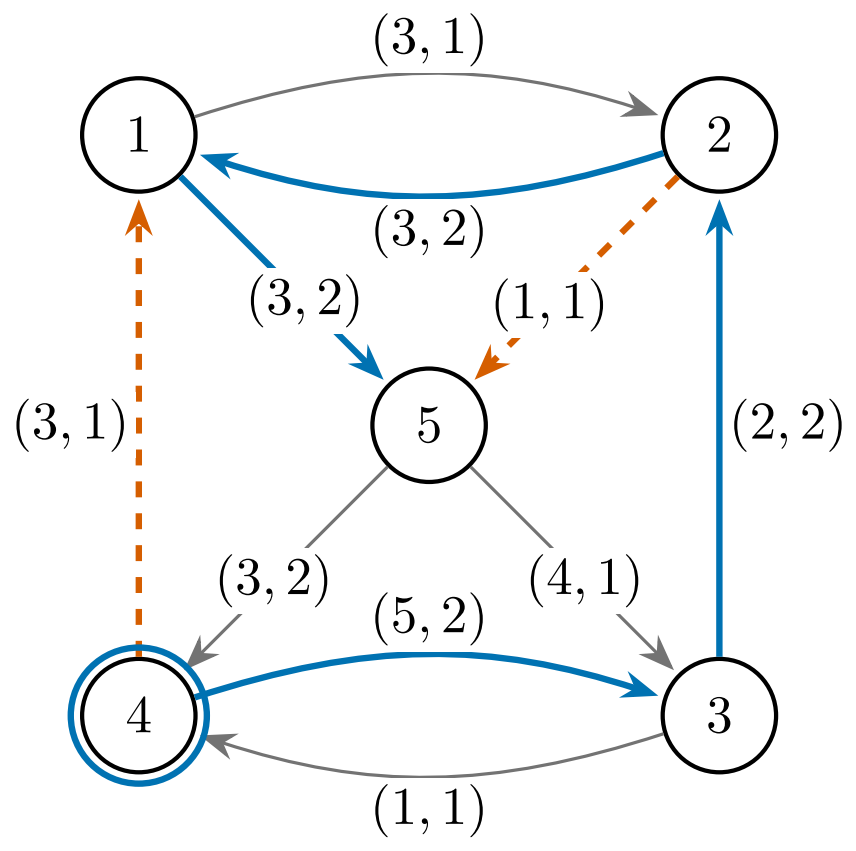
\bigskip
  \begin{subfigure}[t]{0.44\textwidth}\centering\vspace{0pt}
    \begin{tikzpicture}[showcosts, ex/25/.style=interdicted,
                        ex/41/.style=interdicted, ex/43/.style=tree,
                        ex/32/.style=tree, ex/21/.style=tree,
                        ex/15/.style=tree]
      \RunningExample
      \node[root] at (v4) {};
    \end{tikzpicture}
    \caption{$b=2$: arcs~$(2,5)$ and $(4,1)$ are interdicted, and the new shortest paths from vertex~$4$ are shown.}
    \label{FIG_solution_b2}
  \end{subfigure}\hfill
  \begin{subfigure}[t]{0.52\textwidth}\centering\vspace{0pt}
    \begin{tabular}{@{}c|ccccc|c@{}}
      from\,$\backslash$\,to & $1$ & $2$ & $3$ & $4$ & $5$ & sum\\
      \hline
      $1$ & -- & $3$ & $7$ & $6$ & $3$ & $19$\\
      $2$ & $3$ & -- & \chg{10} & \chg{9} & \chg{6} & \chg{28}\\
      $3$ & \chg{5} & $2$ & -- & $1$ & \chg{8} & \chg{16}\\
      $4$ & \chg{10} & \chg{7} & $5$ & -- & \chg{13} & \chg{35}\\
      $5$ & \chg{9} & $6$ & $4$ & $3$ & -- & \chg{22}\\
      \hline
      \multicolumn{6}{@{}r|}{$\psi(\mathcal{G},{\boldsymbol \ell},{\boldsymbol c},2)$} & \chg{120}
    \end{tabular}
    \medskip
    \caption{$b=2$: distances in $\mathcal{G}(\mathcal{V},\mathcal{A}\setminus I)$, $I=\set{(2,5),(4,1)}$.}
    \label{FIG_solution_b2_matrix}
  \end{subfigure}
  \caption{Optimal interdictions of the running example of
    \Cref{FIG_a} for the budgets~$b=1$ (top) and $b=2$ (bottom). Left:
    the interdicted arcs (dashed) and, in blue, the shortest paths from
    the most affected root (circled). For vertex~$2$, compare with
    \Cref{FIG_a_graph}. Right: the distances in the
    interdicted network. The entries, row sums, and APSP values that increase with
    respect to \Cref{FIG_a_matrix} are in red. Each arc is labeled with the pair~$(\ell_a, c_a)$ of its length and interdiction cost.}
  \label{FIG_solution}
\end{figure}

\begin{example}
  \label{ex:running-solution}
  For the running example in~\Cref{FIG_a}, we now set the
  interdiction cost to~$c_a = 2$ for the
  arcs~$(1,5)$, $(2,1)$, $(3,2)$, $(4,3)$,
  and~$(5,4)$, and to $c_a = 1$ for the remaining five
  arcs. In \Cref{FIG_solution}, each arc is labeled with its length and
  its
  cost. Every vertex thus has
  exactly one outgoing and one incoming arc of cost~$2$, and one can
  check that every
  $(s,t)$-cut of the graph has total interdiction cost at least~$3$.
  Hence, no interdiction of cost at most~$2$ can disconnect the graph,
  and Assumptions~\ref{ass:connected} and~\ref{ass:strong-residual} hold
  for $b \in \set{1,2}$. \Cref{FIG_solution} shows the two optimal
  interdictions, which are unique.
  \begin{itemize}
  \item With $b=1$, only the five arcs of cost~$1$ can be interdicted,
    and the resulting values range from $87$ to~$100$. The optimal
    interdiction removes arc~$(2,5)$, which carries three of the four
    shortest paths rooted at vertex~$2$. These paths are rerouted
    through $(2,1)$ and~$(1,5)$, so that the distances from~$2$ to $3$,
    $4$, and~$5$ grow from $5$, $4$, and~$1$ to $10$, $9$, and~$6$. In addition, the distance from~$3$ to~$5$ grows from $3$ to~$7$. Hence,
    $\psi(\mathcal{G},{\boldsymbol \ell},{\boldsymbol c},1) = 100$, against
    $\varphi(\mathcal{G},{\boldsymbol \ell}) = 81$.
  \item With $b=2$, the interdictor can remove either two arcs of
    cost~$1$ or a single arc of cost~$2$. The optimal interdiction
    removes $(2,5)$ and~$(4,1)$, and $\psi(\mathcal{G},{\boldsymbol \ell},{\boldsymbol c},2) = 120$.
    The second arc mostly hurts vertex~$4$, whose paths must now go
    around the square through $(4,3)$, $(3,2)$, and~$(2,1)$: its
    contribution grows from $20$ to~$35$, and its distance to vertex~$5$
    from $6$ to~$13$. The best interdiction of a single arc of cost~$2$,
    namely $(3,2)$, only gives~$93$.
  \end{itemize}
  The example also shows that the effects of interdicted arcs do not
  simply add up: removing $(2,5)$ or $(4,1)$ alone increases the APSP
  value by $19$ and~$14$, respectively, but removing both increases it
  by~$39$. For instance, the distance from~$4$ to~$5$ is $6$ if only
  $(2,5)$ is removed (path~$4 \to 1 \to 5$), $8$ if only $(4,1)$ is
  removed (path~$4 \to 3 \to 2 \to 5$), and $13$ if both are removed.
\end{example}

\subsubsection{Penalty Formulation}

We now introduce the penalty formulation. A value~$M(a,w)\ge0$ is
associated with every arc~$a$ and active root~$w$.

\begin{definition}[Valid penalties]
  \label{def:valid-bigM}
  The values~$M(a,w)$ are \emph{valid} for~$\mathcal K$ if, for every
  feasible interdiction set~$I$ and every requested pair~$(w,t)\in\mathcal K$,
  the shortest $w$--$t$ distance equals $d^I(w,t)$ when each
  arc~$a$ has penalized length~$\ell_a+M(a,w)$ for
  $a\in I$ and $\ell_a$ otherwise.
\end{definition}

Note that validity allows for ties. We call the values strictly valid if no shortest
penalized path of a requested pair uses an interdicted arc.
For a path~$P$, we denote by $\ell(P)\define\sum_{a\in P}\ell_a$ its
length and, for a feasible interdiction set~$I$ and a root~$w$, by
\[
  \ell^M_{I,w}(P)\define\ell(P)+\sum_{a\in P\cap I}M(a,w)
\]
its \emph{penalized length}, i.e., its length when every interdicted
arc~$a$ is lengthened by its penalty~$M(a,w)$.

The penalty formulation is
\begin{subequations}
  \label{eq:apspip-reform}
  \begin{alignat}{3}
    \psi_{\mathcal K}(\mathcal G,{\boldsymbol \ell},{\boldsymbol c},b)
       = \max_{{\boldsymbol x}\in\mathcal X}\quad
       &\makebox[2em][c]{$\displaystyle\min_{\boldsymbol y}$}
       &\quad&\sum_{w\in\mathcal R}\sum_{a\in\mathcal A}
       (\ell_a+M(a,w)x_a)y_{aw}
       &&
       \label{OBJ_apspip_reform}\\*[1.2ex]
       &\makebox[2em][c]{s.t.}
       &&\sum_{a\in\Inedges(t)}y_{aw}
       -\sum_{a\in\Outedges(t)}y_{aw}
       \ge1,
       &\quad&w\in\mathcal R,\ t\in\mathcal T_w,
       \label{C_flow_apspip_reform}\\[1.2ex]
    &&&\sum_{a\in\Inedges(v)}y_{aw}
       -\sum_{a\in\Outedges(v)}y_{aw}
       \ge0,
       &\quad&w\in\mathcal R,\
        v\in\mathcal V\setminus(\mathcal T_w\cup\Set{w}),
       \label{C_zero_apspip_reform}\\[1.2ex]
    &&&y_{aw}\ge0,
       &\quad&w\in\mathcal R,\ a\in\mathcal A.
       \label{C_dom_y_apspip_reform}
  \end{alignat}
\end{subequations}

Formulation~\eqref{eq:apspip-reform} replaces the linking
constraints~\eqref{C_bigM_apspip} by the penalties~$M(a,w)x_a$ in the
follower's arc lengths. Its demand, transit, and nonnegativity
constraints and the leader's feasible set~$\mathcal X$ are as in
Formulation~\eqref{eq:apspip}. The follower now has a linear objective
and a feasible region independent of the interdiction decision, so
linear optimization duality yields a single-level formulation. This
formulation is exact when the penalties are valid.

\section{A Compact Single-Level Formulation}
\label{sec:approaches}
\label{sec:single-level}

Applying LP strong duality to the lower-level problem yields
\begin{subequations}
  \label{eq:single-level}
  \begin{alignat}{3}
    \psi_{\mathcal K}(\mathcal G,{\boldsymbol \ell},{\boldsymbol c},b) =\ 
       &\makebox[2em][c]{$\displaystyle\max_{{\boldsymbol x\in \mathcal X},{\boldsymbol \pi}}$}
       &\quad&\sum_{w\in\mathcal R}\sum_{t\in\mathcal T_w}\pi_{tw}
       &&
       \label{INTER}\\*[1.2ex]
       &\makebox[2em][c]{s.t.}
       &&\pi_{vw}-\pi_{uw}\le \ell_a+M(a,w)x_a,
       &\quad&w\in\mathcal R,\
         a=(u,v)\in\mathcal A\setminus\Outedges(w),
       \label{C1_INTER}\\[1.2ex]
    &&&\pi_{vw}\le \ell_a+M(a,w)x_a,
       &\quad&w\in\mathcal R,\
         a=(w,v)\in\Outedges(w),
       \label{C2_INTER}\\[1.2ex]
    &&&\pi_{vw}\ge0,
       &\quad&w\in\mathcal R,\
         v\in\mathcal V\setminus\Set{w},
       \label{C3_INTER}\\[1.2ex]
    &&&x_a\in\Set{0,1},
       &\quad&a\in\mathcal A.
       \label{C5_INTER}
  \end{alignat}
\end{subequations}

For every active root~$w$, the variable~$\pi_{vw}$ is the dual
variable of the demand or transit constraint of vertex~$v$, and can be
read as a distance label from~$w$ to~$v$. The
objective function~\eqref{INTER} maximizes the sum of the labels of the
requested destinations. The \emph{dual
constraints}~\eqref{C1_INTER} require that the label of the head of
every arc exceed the label of its tail by at most the penalized length
of the arc, and Constraints~\eqref{C2_INTER} impose the same condition on
the arcs leaving the root, whose label is fixed to zero. Hence, by
induction along any path, the label~$\pi_{vw}$ is at most the penalized
length of every $w$--$v$ path.
The budget constraint of the leader is
imposed by the domain~$\mathcal X$ over which the maximization
runs, and
Constraints~\eqref{C3_INTER} and~\eqref{C5_INTER} define the domains of
the variables. The labels are nonnegative because the demand and
transit constraints of the follower are inequalities.
Only requested destination labels appear in the objective. At an
optimal solution, $\pi_{tw}$ equals the penalized shortest-path distance
for every $t\in\mathcal T_w$. Labels of non-destination vertices may
take any values compatible with dual optimality. With valid penalties,
Formulation~\eqref{eq:single-level} is an exact compact formulation of
the MPSP interdiction problem.

The formulation has one dual block per active root and therefore
$O(\card{\mathcal R}m)$ arc constraints. The all-pairs compact model is
recovered when $\mathcal R=\mathcal V$ and
$\mathcal T_w=\mathcal V\setminus\set{w}$.

\section{Valid Penalties}
\label{sec:bigM-analysis}

The characterization of valid penalties depends on which pairs
are requested. In the all-pairs case, for a fixed feasible
interdiction and an active root~$w$, every non-root distance
label appears in the dual objective and therefore equals its
residual distance at optimality. For a general pair set, only the
labels of the requested destinations in~$\mathcal T_w$
are fixed in this way; the remaining non-root labels
may vary among optimal dual solutions.

\subsection{Exact Validity Through the Residual Dual Optimal Face}

Fix a feasible interdiction set~$I$ and an active root~$w\in\mathcal R$. Let $\Pi_w(I)$ be the set of optimal dual potentials
of the residual multiple-pairs subproblem:
\begin{equation}
  \label{eq:dual-optimal-face}
  \Pi_w(I)\define
  \left\{
  {\boldsymbol \pi}_w\in\mathbb R_+^{n-1}:
  \begin{array}{ll}
    \pi_{vw}-\pi_{uw}\le\ell_{(u,v)},
      & (u,v)\in\mathcal A\setminus I,\\[1.5ex]
    \pi_{ww}=0, & \\[1.5ex]
    \displaystyle
    \sum_{t\in\mathcal T_w}\pi_{tw}
    =
    \sum_{t\in\mathcal T_w}d^I(w,t) &
  \end{array}
  \right\},
\end{equation}
i.e., the nonnegative labels from root~$w$ that satisfy the dual
constraints of the residual arcs and whose sum over the requested
destinations equals the sum of the residual distances from~$w$ to them.
By strong duality, $\Pi_w(I)$ is nonempty under
Assumption~\ref{ass:feasibility}.

The following proposition characterizes the valid penalties through these faces.
\begin{proposition}[Exact validity characterization]
  \label{prop:mpsp-bigM-face}
  The penalties~$M(a,w)$ are valid for~$\mathcal K$ if and only if, for
  every feasible interdiction set~$I$ and every active root~$w\in\mathcal R$, there exists
  ${\boldsymbol \pi}_w^I\in\Pi_w(I)$ such that
  \begin{equation}
    \label{eq:mpsp-bigM-face}
    \pi^I_{vw}-\pi^I_{uw}
    \le
    \ell_a+M(a,w),
    \quad
    a=(u,v)\in I.
  \end{equation}
\end{proposition}
\begin{proof}
  Suppose first that such a potential exists. It satisfies the residual
  arc constraints by membership in~$\Pi_w(I)$ and the penalized
  constraints of every interdicted arc by
  \eqref{eq:mpsp-bigM-face}. Hence, it is feasible for the dual of the
  penalized follower and has objective value equal to the residual
  contribution of root~$w$. The penalized follower can only have a
  value no larger than the residual follower because it retains every
  residual path. Strong duality therefore forces equality.

  Conversely, assume that the penalties are valid. For fixed~$I$ and
  $w$, take an optimal dual solution of the penalized follower. Its
  objective equals the residual contribution by validity. Restricted to
  the noninterdicted arcs it is residual-dual feasible, hence belongs to
  $\Pi_w(I)$, and its constraints on the interdicted arcs are exactly
  \eqref{eq:mpsp-bigM-face}.
\end{proof}

Within a single feasible interdiction set, the penalties of different
arcs are coupled: changing a non-destination label can reduce the
required penalty on one arc and increase it on another. Over the whole
family of feasible interdiction sets, however, validity reduces to one
threshold for each arc and root, because a subset of a feasible
interdiction set is feasible. For an interdictable arc~$a=(u,v)$, define
\begin{equation}
  \label{eq:feasible-sets-a}
  \mathcal I(a)\define
  \Set{I\subseteq\mathcal A:
       a\in I,\ \textstyle\sum_{a'\in I}c_{a'}\le b},
\end{equation}
i.e., the family of the feasible interdiction sets that contain~$a$,
and, for every root~$w\in\mathcal R$,
\begin{equation}
  \label{eq:eta-K}
  \eta_{\mathcal K}(a,w)\define
  \max\Set{
    0,\
    \max_{I\in\mathcal I(a)}\max_{t\in\mathcal T_w}
    \bigl(d^I(w,t)-d^I(w,u)-\ell_a-d^I(v,t)\bigr)
  },
\end{equation}
where a term is $-\infty$ if $d^I(w,u)$ or $d^I(v,t)$ is infinite. For the arcs with $c_a>b$, which are never interdicted, we set~$\eta_{\mathcal K}(a,w)\define0$. For a fixed feasible interdiction set~$I$ containing~$a=(u,v)$
and a requested destination~$t$, the sum
$d^I(w,u)+\ell_a+d^I(v,t)$ is the length obtained by following a
shortest surviving path from~$w$ to~$u$, traversing~$a$, and following
a shortest surviving path from~$v$ to~$t$, whenever both surviving
paths exist. The amount by which this length falls short of~$d^I(w,t)$
determines the penalty needed to prevent this route from being shorter
than the residual shortest path. The
value~$\eta_{\mathcal K}(a,w)$ takes the largest required
penalty over all requested destinations and feasible interdiction sets
containing~$a$.

The following proposition shows that $\eta_{\mathcal K}$ is the componentwise smallest valid penalty vector, for every pair set.
\begin{proposition}[Componentwise smallest penalties]
  \label{prop:mpsp-smallest}
  Under Assumption~\ref{ass:feasibility}, the penalties are valid if
  and only if $M(a,w)\geq\eta_{\mathcal K}(a,w)$ for every
  interdictable arc~$a$ and every root~$w\in\mathcal R$. Hence,
  $\eta_{\mathcal K}$ is the componentwise smallest valid penalty
  vector.
\end{proposition}
\begin{proof}
  Assume that the penalties are valid, and let $I\in\mathcal I(a)$ and
  $t\in\mathcal T_w$ with $d^I(w,u)$ and $d^I(v,t)$ finite. A shortest
  residual $w$--$u$ path, arc~$a$, and a shortest residual $v$--$t$ path
  form a $w$--$t$ walk whose only interdicted arc is~$a$. Since the
  lengths are nonnegative, the penalized $w$--$t$ distance is at most
  the penalized length of this walk, and it equals $d^I(w,t)$ by
  validity. Hence,
  $d^I(w,u)+\ell_a+M(a,w)+d^I(v,t)\geq d^I(w,t)$, and
  $M(a,w)\geq\eta_{\mathcal K}(a,w)$ because $M(a,w)\geq0$.

  Conversely, assume $M(a,w)\geq\eta_{\mathcal K}(a,w)$, and fix a
  feasible interdiction set~$I$, a root~$w$, and $t\in\mathcal T_w$. A
  shortest residual $w$--$t$ path has penalized length~$d^I(w,t)$, so it
  suffices to show that $\ell^M_{I,w}(P)\geq d^I(w,t)$ for every simple
  $w$--$t$ path~$P$. We proceed by induction on $k=\card{P\cap I}$. For
  $k=0$, $P$ is residual and the claim holds. For $k\geq1$, let
  $a_1,\ldots,a_k$ be the interdicted arcs of~$P$ in path order, let
  $a_k=(u,v)$, and let
  $J\define I\setminus\set{a_1,\ldots,a_{k-1}}$, i.e., the set
  obtained from~$I$ by restoring the first $k-1$ interdicted arcs
  of~$P$, while all other interdictions stay in place. Since the costs are
  positive, $J\in\mathcal I(a_k)$. The subpath of~$P$ from~$w$ to~$u$
  avoids~$J$, and the subpath from~$v$ to~$t$ avoids~$I\supseteq J$. Hence, their lengths are at least $d^J(w,u)$ and $d^J(v,t)$,
  respectively. By~\eqref{eq:eta-K} applied to~$J$,
  $M(a_k,w)\geq d^J(w,t)-d^J(w,u)-\ell_{a_k}-d^J(v,t)$, and therefore
  \[
    \ell^M_{I,w}(P)\geq d^J(w,t)+\sum_{i=1}^{k-1}M(a_i,w).
  \]
  By Assumption~\ref{ass:feasibility}, there is a shortest $w$--$t$
  path~$Q$ in the interdicted network of~$J$. Its arcs in~$I$ belong to
  $\set{a_1,\ldots,a_{k-1}}$, so that, by induction and since the
  penalties are nonnegative,
  \[
    d^I(w,t)\leq\ell^M_{I,w}(Q)
    =d^J(w,t)+\sum_{a'\in Q\cap I}M(a',w)
    \leq d^J(w,t)+\sum_{i=1}^{k-1}M(a_i,w)
    \leq\ell^M_{I,w}(P).\qedhere
  \]
\end{proof}

The proof removes the coupling of a single interdiction set by passing
to the smaller set~$J$, in which the prefix of the path is residual. Only distances to requested destinations are used. The set of valid
penalties is thus closed under componentwise minima. The
characterization also gives a way to construct valid penalties: any
upper bound on $\eta_{\mathcal K}(a,w)$ is a valid penalty, and upper
bounds on residual distances provide such values without enumerating
the feasible interdiction sets. We compute $\eta_{\mathcal K}$ itself by
enumeration on small instances and use it only as a theoretical
reference.

\subsection{The All-Pairs Penalties as a Special Case}

For $\mathcal K=\mathcal K^{\mathrm{AP}}$, every non-root vertex is a
destination. Any residual-dual feasible label satisfies
$\pi_{vw}\le d^I(w,v)$. Equality of the sums in
\eqref{eq:dual-optimal-face} then forces equality componentwise. Thus,
\[
  \Pi_w(I)=\Set{(d^I(w,v))_{v\neq w}}.
\]

Under Assumption~\ref{ass:strong-residual}, define
\begin{equation}
  \label{eq:rooted-gap}
  \eta(a,w)\define
  \max\Set{
    0,\
    \max_{I\in\mathcal I(a)}
    \bigl(d^I(w,v)-d^I(w,u)-\ell_a\bigr)
  }.
\end{equation}
The value~$\eta(a,w)$ is the largest amount by which a shortest residual
$w$--$u$ path followed by~$a$ falls short of the residual $w$--$v$
distance, over the sets in~$\mathcal I(a)$.
For $\mathcal K=\mathcal K^{\mathrm{AP}}$ and $w\neq v$, the residual
triangle inequality $d^I(w,t)\leq d^I(w,v)+d^I(v,t)$ shows that the
inner maximum in~\eqref{eq:eta-K} is attained at $t=v$. For $w=v$,
every term in~\eqref{eq:eta-K} and in~\eqref{eq:rooted-gap} equals
$-d^I(w,u)-\ell_a\le0$, since $d^I(v,v)=0$, so that both values are
zero. Hence, $\eta_{\mathcal K^{\mathrm{AP}}}=\eta$, where $\eta(a,w)\define0$
for the arcs with $c_a>b$ as well.

Proposition~\ref{prop:mpsp-smallest} then gives the following characterization of the all-pairs penalties.
\begin{corollary}[All-pairs componentwise smallest penalties]
  \label{prop:rooted-bigM}
  For $\mathcal K=\mathcal K^{\mathrm{AP}}$, the penalties are valid if
  and only if
  \[
    M(a,w)\ge\eta(a,w)
  \]
  for every interdictable arc~$a$ and every root~$w$.
\end{corollary}

For root-independent values,
\begin{equation}
  \label{eq:robust-detour}
  \eta(a)\define
  \max_{w\in\mathcal V}\eta(a,w)
  =
  \eta(a,u)
  =
  \max\Set{
    0,\
    \max_{I\in\mathcal I(a)}d^I(u,v)-\ell_a
  }.
\end{equation}
The value~$\eta(a)$ is the largest amount by which~$\ell_a$ falls short
of the residual $u$--$v$ distance, over the sets in~$\mathcal I(a)$.
The second equality holds because the residual triangle inequality
$d^I(w,v)\leq d^I(w,u)+d^I(u,v)$ gives
$d^I(w,v)-d^I(w,u)\leq d^I(u,v)$, with equality for $w=u$.
Taking the maximum over the roots gives the root-independent counterpart of Corollary~\ref{prop:rooted-bigM}.
\begin{corollary}[All-pairs root-independent penalties]
  \label{prop:robust-arc-bigM}
  For $\mathcal K=\mathcal K^{\mathrm{AP}}$, root-independent
  penalties~$M(a,w)=M(a)$ are valid if and only if
  $M(a)\ge\eta(a)$ for every interdictable arc~$a$.
\end{corollary}

Hence, $M(a)=\eta(a)$ is the componentwise smallest root-independent
all-pairs value.

\begin{example}
  \label{ex:running-exact}
  We compute the smallest valid all-pairs values on the running example
  of \Cref{FIG_a}. \Cref{FIG_bigM_exact} illustrates the two
  quantities for the arc~$(1,2)$, and \Cref{tab:running-exact} reports
  them for all arcs. For $b=2$, the root-dependent values sum to~$89$,
  against $5\cdot53=265$ for the root-independent ones.
\end{example}

The smallest valid values are a theoretical reference only: we compute
them by enumeration on small instances, and such an enumeration would
also solve the problem, so they are not used as rules.

\begin{figure}[htb]
  \centering
  \captionsetup[subfigure]{font=wiley}
  \begin{subfigure}[t]{0.45\textwidth}\centering\vspace{0pt}
    \begin{tikzpicture}[showcosts, ex/12/.style=interdicted, ex/15/.style=detour,
                        ex/53/.style=detour, ex/32/.style=detour]
      \RunningExample
    \end{tikzpicture}
    \caption{The shortest detour of $(1,2)$ has length~$9$: for $b=1$,
      $\eta((1,2))=\max\set{0,9-\ell_{(1,2)}}=6$, so that
      Corollary~\ref{prop:robust-arc-bigM} requires
      $M((1,2))\geq6$.}
    \label{FIG_bigM_detour}
  \end{subfigure}
\hfill
  \begin{subfigure}[t]{0.45\textwidth}\centering\vspace{0pt}
    \begin{tikzpicture}[showcosts, ex/12/.style=interdicted, ex/41/.style=tree,
                        ex/43/.style=tree, ex/32/.style=tree]
      \RunningExample
      \node[root] at (v4) {};
    \end{tikzpicture}
    \caption{Root~$4$ with $b=2$: without $(1,2)$, vertex~$4$ reaches
      vertex~$1$ at distance~$3$ and vertex~$2$ at distance~$7$, so
      $\eta((1,2),4)=\max\set{0,7-3-\ell_{(1,2)}}=1$ and
      Corollary~\ref{prop:rooted-bigM} requires only
      $M((1,2),4)\geq1$, against $M((1,2),1)\geq10$ for the root~$1$.}
    \label{FIG_bigM_rooted}
  \end{subfigure}
  \caption{Smallest valid penalties for the interdicted arc~$(1,2)$
    (dashed) in the running example: the worst-case excess detour penalty
    of Corollary~\ref{prop:robust-arc-bigM} and, for a given root
    (circled), the residual-distance penalty of
    Corollary~\ref{prop:rooted-bigM}.}
  \label{FIG_bigM_exact}
\end{figure}

\begin{table}[htb]
  \centering
  \caption{Smallest valid penalties of the running example:
    $\eta(a)$ for $b=1$ and $b=2$, and $\eta(a,w)$ for $b=2$ and every
    root~$w$ (a dash
    means that the arc cannot be interdicted).}
  \label{tab:running-exact}
  \small
  \begin{tabular}{@{}crr r rrrrr@{}}
    \toprule
    & \multicolumn{2}{c}{$\eta(a)$} & & \multicolumn{5}{c}{$\eta(a,w)$, $b=2$}\\
    \cmidrule{2-3}\cmidrule{5-9}
    Arc~$a=(u,v)$ & $b=1$ & $b=2$ & & $w=1$ & $2$ & $3$ & $4$ & $5$\\
    \midrule
    $(1,2)$ &  6 & 10 & & 10 &  0 &  0 &  1 &  1\\
    $(1,5)$ & -- &  1 & &  1 &  0 &  0 &  1 &  0\\
    $(2,1)$ & -- &  4 & &  0 &  4 &  0 &  0 &  0\\
    $(2,5)$ &  5 &  5 & &  0 &  5 &  5 &  5 &  0\\
    $(3,2)$ & -- &  5 & &  0 &  0 &  5 &  0 &  3\\
    $(3,4)$ &  5 & 10 & &  0 &  0 & 10 &  0 &  0\\
    $(4,1)$ &  7 &  7 & &  0 &  0 &  1 &  7 &  7\\
    $(4,3)$ & -- &  5 & &  0 &  0 &  0 &  5 &  0\\
    $(5,3)$ &  4 &  4 & &  4 &  4 &  0 &  0 &  4\\
    $(5,4)$ & -- &  2 & &  2 &  2 &  0 &  0 &  2\\
    \bottomrule
  \end{tabular}
\end{table}

For a general pair set and under
Assumption~\ref{ass:strong-residual}, the same triangle inequality
gives $\eta_{\mathcal K}(a,w)\leq\eta(a,w)$. Hence, the all-pairs values~$\eta(a,w)$ remain valid penalties for every pair set, and the
destination-aware values~$\eta_{\mathcal K}(a,w)$ can be strictly
smaller, as the following example shows.

\begin{example}
  \label{ex:running-pairs}
  On the running example of \Cref{FIG_a}, let $b=2$ and
  $\mathcal K=\set{(1,4),(1,5)}$, so that $\mathcal R=\set{1}$ and
  $\mathcal T_1=\set{4,5}$. For $a=(1,2)$ and $w=1$, the all-pairs value
  is $\eta(a,1)=10$: removing $a$ together with $(5,3)$ makes the
  residual distance from~$1$ to~$2$ equal to~$13=\ell_a+10$. By contrast,
  $\eta_{\mathcal K}(a,1)=0$. Indeed, every feasible interdiction set
  containing~$a$ leaves the arcs~$(1,5)$ and $(5,4)$, whose cost is~$2$,
  in place, so that the requested distances remain $d(1,5)=3$ and
  $d(1,4)=6$. In the interdiction $\set{(1,2),(5,3)}$, the labels~$(\pi_{11},\pi_{21},\pi_{31},\pi_{41},\pi_{51})=(0,3,5,6,3)$ satisfy
  the dual constraints of every arc with zero penalties, those of the
  two interdicted arcs included, and the requested labels sum to the
  residual value~$9$,
  although $\pi_{21}=3$ differs from $d^I(1,2)=13$: the label of a vertex
  that is not requested need not equal its residual distance. Over all
  arcs, $\eta_{\mathcal K}(a,1)$ equals~$1$ on $(1,5)$, $2$ on $(5,4)$,
  and~$0$ on the other eight arcs.
\end{example}

\subsection{Poly-Time Computable Penalty Rules}
\label{sec:bigM-bounds}

We now derive four poly-time computable penalty rules from upper bounds on the
residual distances. Two of them are root-independent, and the other two
exploit the root.

\subsubsection{Arc-length bounds}

For an arc~$a=(u,v)$, define
\[
  \mathcal{A}(a)\define
  \mathcal{A}\setminus\bigl(\Set{a}\cup\Inedges(u)\cup\Outedges(v)\bigr),
\]
i.e., the arcs other than~$a$ that neither enter~$u$ nor leave~$v$.
Only these arcs can belong to a simple $u$--$v$ path avoiding~$a$.

The following corollary bounds the residual $u$--$v$ distance by the sum of the largest lengths of these arcs.
\begin{corollary}
  \label{cor:bigM-Lbar}
  Let Assumptions~\ref{ass:connected} and~\ref{ass:strong-residual} be
  satisfied and, for every arc~$a$, let
  $m_a\define|\mathcal{A}(a)|$ and let
  $a_1(a),\ldots,a_{m_a}(a)$ be an ordering of~$\mathcal{A}(a)$ satisfying
  $\ell_{a_1(a)}\geq\cdots\geq\ell_{a_{m_a}(a)}$. Define the
  arc-length bound
  \begin{equation}
    \label{eq:sigma}
    \sigma(a)\define
      \sum_{j=1}^{\min\Set{n-1,m_a}}\ell_{a_j(a)},
  \end{equation}
  i.e., the total length of the $\min\set{n-1,m_a}$ longest arcs
  of~$\mathcal{A}(a)$. Then, $d^I(u,v)\leq\sigma(a)$ for every interdictable arc~$a=(u,v)$
  and every $I\in\mathcal I(a)$.
\end{corollary}

The bound of Corollary~\ref{cor:bigM-Lbar} only requires sorting the arc lengths once, in $O(m\log m)$ time; then $\sigma(a)$ is obtained for every arc $a=(u,v)$ by masking out $\{a\}\cup\Inedges(u)\cup\Outedges(v)$ and summing the $\min\{n-1,m_a\}$ longest remaining arcs in the global order, in $O(n+\deg^-(u)+\deg^+(v))$ time per arc.

For arc~$(2,1)$ in the running example, $\sigma((2,1))=15$.

The root-dependent bound also accounts for the distance already paid
from the root to the tail. A residual shortest $w$--$v$ path is simple,
neither enters~$w$ nor leaves~$v$, and satisfies
$d^I(w,u)\geq d(w,u)$. For $w\neq v$, define
\[
  \mathcal{A}'(w,v)\define\mathcal{A}\setminus
  \bigl(\Inedges(w)\cup\Outedges(v)\bigr),
  \quad r_{w,v}\define|\mathcal{A}'(w,v)|,
\]
let $a'_1(w,v),\ldots,a'_{r_{w,v}}(w,v)$ be an ordering of
$\mathcal{A}'(w,v)$ satisfying
$\ell_{a'_1(w,v)}\geq\cdots\geq\ell_{a'_{r_{w,v}}(w,v)}$, and define
\[
  \sigma(w,v)\define
  \sum_{j=1}^{\min\Set{n-1,r_{w,v}}}\ell_{a'_j(w,v)},
\]
i.e., the total length of the $\min\set{n-1,r_{w,v}}$ longest arcs
that neither enter~$w$ nor leave~$v$. In the same way, $d^I(w,v)\leq\sigma(w,v)$ for every feasible
interdiction set~$I$.

\subsubsection{Bounds from Surviving Paths}
\label{sec:bigM-flow-bounds}

The arc-length bounds can be tightened with families of paths of which
at least one survives every feasible interdiction. Let
$c^{\min}\define\min_{a\in\mathcal{A}}c_a$. For an interdictable arc~$a=(u,v)$, let $h_a$ be the number of arcs~$a'\neq a$ with
$c_{a'}\leq b-c_a$, and define
\[
  k_a\define\min\Set{\left\lfloor\frac{b-c_a}{c^{\min}}\right\rfloor,\ h_a},
\]
i.e., an upper bound on the number of arcs that a feasible
interdiction set can contain besides~$a$, and let $\mathcal{G}^{F}_a$ be the capacitated graph with arc set~$\mathcal{A}\setminus\set{a}$ and capacities
\[
  q^a_{a'}\define
  \begin{cases}
    1, & c_{a'}\leq b-c_a,\\
    k_a+1, & c_{a'}>b-c_a.
  \end{cases}
\]
The first capacity applies to arcs that can be interdicted together with
$a$. The second applies to arcs that cannot be interdicted together with
$a$ within the budget. Hence, in~$\mathcal{G}^{F}_a$, every arc that can
be interdicted together with~$a$ carries at most one unit of flow.

On $\mathcal{G}^{F}_a$, solve the minimum-cost flow problem
\begin{equation}
  \label{eq:surviving-flow-a}
  \begin{aligned}
    \min\quad & \sum_{a'\in\mathcal{A}\setminus\Set{a}}\ell_{a'} f_{a'}\\
    \text{s.t.}\quad
      &\sum_{a'\in\Outedges(i)\setminus\Set{a}}f_{a'}
       -\sum_{a'\in\Inedges(i)\setminus\Set{a}}f_{a'}
       =\begin{cases}
          k_a+1, & i=u,\\
          -(k_a+1), & i=v,\\
          0, & i\in\mathcal{V}\setminus\Set{u,v},
        \end{cases}
        && i\in\mathcal{V},\\
      & 0\leq f_{a'}\leq q^a_{a'},\quad f_{a'}\in\mathbb{Z},
        && a'\in\mathcal{A}\setminus\Set{a}.
  \end{aligned}
\end{equation}
Let $\mathcal{P}_a$ be any decomposition of an optimal integral flow into
$k_a+1$ $u$--$v$ paths, which need not be distinct, and define
\begin{equation*}
  \mu(a)\define
  \max_{P\in\mathcal{P}_a}\ell(P),
\end{equation*}
i.e., the length of the longest of $k_a+1$ $u$--$v$ paths. The flow
objective selects a family of paths of minimum total length, and the
longest path of the family bounds the length of the path that
survives. We show
in Proposition~\ref{lem:disjoint-paths} below that at least one of them
survives every feasible interdiction containing~$a$. If the
flow problem is infeasible, set~$\mu(a)=+\infty$.

\begin{example}[The surviving-path flow]
  Consider the running example with budget~$b=2$ and let $a=(1,2)$,
  whose length and interdiction cost are $\ell_a=3$ and $c_a=1$.
  Since $c^{\min}=1$, we have
  $k_a=\lfloor(2-1)/1\rfloor=1$, so the auxiliary flow must send
  $k_a+1=2$ units from vertex~$1$ to vertex~$2$ in
  $\mathcal{G}^{F}_{(1,2)}$. The arcs with cost at most
  $b-c_a=1$ receive capacity~$1$. These are
  $(2,5)$, $(3,4)$, $(4,1)$, and $(5,3)$. Every other arc receives
  capacity~$k_a+1=2$. In particular, the two units may share the
  arcs~$(1,5)$ and $(3,2)$, both of which have interdiction cost~$2$
  and therefore cannot be interdicted together with~$a$.

  An optimal integral flow decomposes into the two paths
  \[
    P^1: 1\to5\to3\to2,
    \quad \ell(P^1)=3+4+2=9,
  \]
  and
  \[
    P^2: 1\to5\to4\to3\to2,
    \quad \ell(P^2)=3+3+5+2=13.
  \]
  They share only arcs of capacity~$2$, while every arc of capacity~$1$
  belongs to at most one of them. Any feasible interdiction containing
  $a$ can contain at most one additional arc, and that additional arc
  has capacity~$1$. It can therefore destroy at most one of the two
  paths, so at least one path survives. Thus,
  $\mu((1,2))=\max\set{9,13}=13$.
\end{example}

\begin{figure}[htb]
  \centering
  \begin{tikzpicture}[scale=0.9, transform shape]
    \SurvivingFlowGraph
  \end{tikzpicture}
  \caption{The surviving-path flow for $a=(1,2)$ and $b=2$ in the
    capacitated graph~$\mathcal{G}^{F}_{(1,2)}$. The arc~$a=(1,2)$ is
    absent. The remaining arcs are labeled with $(\ell_{a'},q^a_{a'})$. The
    two flow paths~$1\to5\to3\to2$ and $1\to5\to4\to3\to2$ are shown
    in blue and orange, respectively, and their shared arcs in green.}
  \label{FIG_surviving-flow}
\end{figure}

The root-dependent construction is analogous. For $w\neq v$, let
\[
  k_0\define\min\Set{\left\lfloor\frac{b}{c^{\min}}\right\rfloor,\ h_0},
\]
where $h_0$ is the number of arcs~$a$ with $c_a\leq b$. Thus, $k_0$
is an upper bound on the number of arcs of a feasible interdiction set.
Let $\mathcal{G}^{F}_0$ be the capacitated graph with arc set~$\mathcal{A}$ and capacities
\[
  q^0_a\define
  \begin{cases}
    1, & c_a\leq b,\\
    k_0+1, & c_a>b,
  \end{cases}
\]
so that every arc that can be interdicted carries at most one unit of
flow. Solve the analogous minimum-cost flow problem on $\mathcal{G}^{F}_0$,
with source~$w$, sink~$v$, and $k_0+1$ units of flow, using the capacities~$q^0$ and all arcs in~$\mathcal{A}$.
For an optimal integral flow, let $\mathcal{P}_{w,v}$ be any decomposition
into $k_0+1$ paths and define
\[
  \mu(w,v)\define\max_{P\in\mathcal{P}_{w,v}}\ell(P),
\]
and set~$\mu(w,v)=+\infty$ if the flow problem is infeasible. As
for~$\mu(a)$, the value~$\mu(w,v)$ is the length of the longest of
$k_0+1$ $w$--$v$ paths.

The following proposition shows that at least one of these paths
survives every feasible interdiction, so that the flows bound the
residual distances.
\begin{proposition}
  \label{lem:disjoint-paths}
  For every interdictable arc~$a=(u,v)$ and every
  $I\in\mathcal I(a)$, $d^I(u,v)\leq\mu(a)$. For every ordered pair of
  distinct vertices~$w$ and~$v$ and every feasible interdiction
  set~$I$, $d^I(w,v)\leq\mu(w,v)$.
\end{proposition}
\begin{proof}
  The bounds hold trivially when the flows are infeasible. A feasible
  interdiction set~$I\ni a$ contains, besides~$a$, at most
  $\lfloor(b-c_a)/c^{\min}\rfloor$ arcs, all with cost at most $b-c_a$,
  hence at most $k_a$ arcs,
  and each of them has capacity~$1$, so that it lies on at most one of
  the paths that define~$\mu(a)$. Hence, one of these paths avoids~$I$,
  and $d^I(u,v)\leq\mu(a)$. In the same way, a feasible~$I$ contains
  at most $k_0$ arcs, each with capacity~$1$ in~$\mathcal{G}^{F}_0$, so
  that one of the paths that define~$\mu(w,v)$ avoids~$I$ and
  $d^I(w,v)\leq\mu(w,v)$.
\end{proof}

The proposition holds for every minimum-cost flow and every path
decomposition. With unit costs and $b=1$, the flow has one path and
$\mu(a)=d^{\set{a}}(u,v)$.
Assumption~\ref{ass:strong-residual} guarantees that these flows exist
with unit interdiction costs. With general costs, $\mu$ may be
infinite.

\subsubsection{The Four Poly-Time Computable Rules}
\label{sec:bigM-summary}

For every interdictable arc~$a=(u,v)$ and every root~$w\neq v$, we
define the four poly-time computable penalty rules as
\begin{align}
  M_{\texttt{top}}(a)
    &\define\max\Set{0,\ \sigma(a)-\ell_a},
    \label{eq:M-top}\\
  M_{\texttt{flow}}(a)
    &\define\max\Set{0,\ \min\Set{\mu(a),\sigma(a)}-\ell_a},
    \label{eq:M-flow}\\
  M_{\texttt{rtop}}(a,w)
    &\define\max\Set{0,\ \min\Set{\sigma(a),\ \sigma(w,v)-d(w,u)}-\ell_a},
    \label{eq:M-rtop}\\
  M_{\texttt{rflow}}(a,w)
    &\define\max\Set{0,\ \min\Set{\mu(a),\ \sigma(a),\
      \mu(w,v)-d(w,u),\ \sigma(w,v)-d(w,u)}-\ell_a}.
    \label{eq:M-rflow}
\end{align}
That is, each rule is the amount, if positive, by which the smallest of
its upper bounds exceeds~$\ell_a$: $M_{\texttt{top}}$ uses the
arc-length bound~$\sigma(a)$ on the residual $u$--$v$ distance,
$M_{\texttt{flow}}$ adds the surviving-path bound~$\mu(a)$, and
$M_{\texttt{rtop}}$ and $M_{\texttt{rflow}}$ add the corresponding
bounds on the residual $w$--$v$ distance, reduced by the distance~$d(w,u)$ from the root to the tail of~$a$.
Since $\sigma$ is finite, an infinite value of~$\mu$ does not affect
these minima. The rules $M_{\texttt{top}}$ and $M_{\texttt{flow}}$ are
root-independent, i.e., $M(a,w)=M_{\texttt{top}}(a)$ or
$M(a,w)=M_{\texttt{flow}}(a)$ for every root~$w$. We set all four
values to zero for the arcs with $c_a>b$, which are never interdicted,
and the root-dependent values to zero for $w=v$. By construction,
\[
  M_{\texttt{rflow}}\leq M_{\texttt{flow}}\leq M_{\texttt{top}}
  \quad\text{and}\quad
  M_{\texttt{rflow}}\leq M_{\texttt{rtop}}\leq M_{\texttt{top}}
\]
componentwise. The rules $M_{\texttt{flow}}$ and $M_{\texttt{rtop}}$
are not comparable in general: in the running example with $b=2$
(\Cref{tab:running-bigM}), $M_{\texttt{flow}}((3,4))=10$, while
$M_{\texttt{rtop}}((3,4),1)=5$ and $M_{\texttt{rtop}}((3,4),3)=11$.

The following proposition states that the four rules are valid for
every pair set. It only requires the requested pairs to remain
connected, so that the auxiliary distances~$d^I(u,v)$ and $d^I(w,v)$
used to motivate the bounds may be infinite.
{
Because $M_{\texttt{rflow}}$ is componentwise no larger than the
other three rules, proving validity for $M_{\texttt{rflow}}$ proves
validity for all four rules.}

\begin{proposition}
  \label{prop:polynomial-rules}
  Under Assumption~\ref{ass:feasibility}, each of the four rules
  satisfies $M(a,w)\geq\eta_{\mathcal K}(a,w)$ for every interdictable
  arc~$a$ and every root~$w\in\mathcal R$, and is therefore valid
  for~$\mathcal K$.
\end{proposition}
\begin{proof}
  Since $M_{\texttt{rflow}}$ is the smallest of the four rules, it
  suffices to show that $\eta_{\mathcal K}(a,w)\leq
  M_{\texttt{rflow}}(a,w)$, i.e., that
  \[
    g\define d^I(w,t)-d^I(w,u)-\ell_a-d^I(v,t)
  \]
  is at most each of $\sigma(a)-\ell_a$, $\mu(a)-\ell_a$,
  $\sigma(w,v)-d(w,u)-\ell_a$, and $\mu(w,v)-d(w,u)-\ell_a$, for every
  $I\in\mathcal I(a)$ and $t\in\mathcal T_w$ with $d^I(w,u)$ and
  $d^I(v,t)$ finite. For $w=v$, we have $g=-d^I(v,u)-\ell_a\le0$, since
  $d^I(v,v)=0$, so that $\eta_{\mathcal K}(a,v)=0$. Let $w\neq v$. By
  Assumption~\ref{ass:feasibility}, there is a shortest $w$--$t$
  path~$P$ in the interdicted network, which is simple and avoids~$a$. Its subpaths are shortest residual paths.

  \emph{Bound~$\sigma(a)$.} If $u$ and~$v$ lie on~$P$ in this order, the
  prefix of~$P$ up to~$u$ and its suffix from~$v$ have lengths~$d^I(w,u)$ and $d^I(v,t)$, so that $g$ is the length of the
  $u$--$v$ subpath of~$P$ minus~$\ell_a$. If $v$ precedes~$u$, then
  $d^I(w,t)\le d^I(w,u)+d^I(v,t)$ and $g\le-\ell_a$. If only~$u$ lies
  on~$P$, then $g$ is at most the length of the $u$--$t$ subpath
  minus~$\ell_a$. If only~$v$ lies on~$P$, at most the length of the
  $w$--$v$ subpath minus~$\ell_a$. If neither does, at most
  $\ell(P)-\ell_a$. In each case, the subpath is simple, avoids~$a$,
  neither enters~$u$ nor leaves~$v$, and thus has at most $n-1$ arcs
  of~$\mathcal A(a)$, so that its length is at most~$\sigma(a)$.

  \emph{Bound~$\sigma(w,v)$.} If $v$ lies on~$P$, then
  $d^I(w,t)=d^I(w,v)+d^I(v,t)$ and the $w$--$v$ subpath of~$P$ neither
  enters~$w$ nor leaves~$v$. Otherwise, $P$ itself neither enters~$w$
  nor leaves~$v$ and $d^I(w,t)\le\ell(P)$. In both cases, since
  $d^I(w,u)\geq d(w,u)$ and $d^I(v,t)\ge0$, we have
  $g\le\sigma(w,v)-d(w,u)-\ell_a$.

  \emph{Bounds~$\mu(a)$ and $\mu(w,v)$.} If $\mu(a)$ is finite, by
  Proposition~\ref{lem:disjoint-paths} one of the paths that define it
  is a residual $u$--$v$ path of length at most~$\mu(a)$, so that
  $d^I(w,t)\le d^I(w,u)+\mu(a)+d^I(v,t)$ and $g\le\mu(a)-\ell_a$. In
  the same way, if $\mu(w,v)$ is finite, then
  $d^I(w,t)\le\mu(w,v)+d^I(v,t)$ and $g\le\mu(w,v)-d(w,u)-\ell_a$. An
  infinite value imposes no condition.

  Hence, $g$ is at most the minimum in~\eqref{eq:M-rflow} minus~$\ell_a$,
  and maximizing over~$I$ and~$t$ gives $\eta_{\mathcal K}(a,w)\le
  M_{\texttt{rflow}}(a,w)$. Validity follows from
  Proposition~\ref{prop:mpsp-smallest}.
\end{proof}

For $\mathcal K=\mathcal K^{\mathrm{AP}}$,
Assumption~\ref{ass:feasibility} coincides with
Assumption~\ref{ass:strong-residual} and $\eta_{\mathcal K}=\eta$, so
that the four rules dominate the all-pairs values~$\eta(a,w)$ of
Corollary~\ref{prop:rooted-bigM}. Taking the minimum of several bounds
preserves validity because each bound dominates
$\eta_{\mathcal K}(a,w)$, and by Proposition~\ref{prop:mpsp-smallest}
the valid penalties are exactly those above~$\eta_{\mathcal K}$. The
four rules ignore the requested destinations, which
$\eta_{\mathcal K}$ exploits.

\begin{figure}[htb]
  \centering
  \captionsetup[subfigure]{font=wiley}
  \begin{subfigure}[t]{\PenaltyPanelWidth}\centering\vspace{0pt}
    \begin{tikzpicture}[scale=0.8, transform shape, showcosts,
                        ex/21/.style=interdicted]
      \RunningExample
    \end{tikzpicture}
    \caption{The arc-length rule gives
      $M_{\texttt{top}}((2,1))=15-\ell_{(2,1)}=12$.}
    \label{FIG_bigM_top}
  \end{subfigure}\hfill
  \begin{subfigure}[t]{\PenaltyPanelWidth}\centering\vspace{0pt}
    \begin{tikzpicture}[scale=0.8, transform shape, showcosts,
                        ex/21/.style=interdicted, ex/25/.style=detour,
                        ex/54/.style=detour, ex/41/.style=detour]
      \RunningExample
    \end{tikzpicture}
    \caption{For $b=2$, no arc can be interdicted together with $(2,1)$,
      and the shortest path suffices: $\mu((2,1))=7$ and
      $M_{\texttt{flow}}((2,1))=7-\ell_{(2,1)}=4$ by
      Proposition~\ref{lem:disjoint-paths}.}
    \label{FIG_bigM_flow1}
  \end{subfigure}\PenaltyPanelBreak
  \begin{subfigure}[t]{\PenaltyPanelWidth}\centering\vspace{0pt}
    \begin{tikzpicture}[scale=0.8, transform shape, showcosts,
                        ex/12/.style=interdicted, ex/15/.style=tree,
                        ex/32/.style=tree, ex/53/.style=detour,
                        ex/54/.style=longest, ex/43/.style=longest]
      \RunningExample
    \end{tikzpicture}
    \caption{For $b=2$, two $1$--$2$ paths share only arcs of cost~$2$. The shared arcs are blue and the path-specific arcs are green and
      orange. They give $\mu((1,2))=13$ and
      $M_{\texttt{flow}}((1,2))=13-\ell_{(1,2)}=10$ by
      Proposition~\ref{lem:disjoint-paths}.}
    \label{FIG_bigM_flow}
  \end{subfigure}
  \caption{Penalty bounds for interdicted arcs (dashed) in the running
  example of \Cref{FIG_a}. Arcs are labeled (length, cost). The first
  panel illustrates $M_{\texttt{top}}$, while the other panels illustrate
  $M_{\texttt{flow}}$.}
  \label{FIG_bigM}
\end{figure}

\begin{example}
  \label{ex:running-bigM}
  In the running example with $b=1$, only the arcs of cost~$1$ can be
  interdicted, one at a time, while
  $M_{\texttt{top}}(a)=M_{\texttt{flow}}(a)=0$ for the arcs of cost~$2$.
  For arc~$(1,2)$, Proposition~\ref{lem:disjoint-paths} uses a single
  path, the shortest $1$--$2$ path without $(1,2)$, namely
  $1 \to 5 \to 3 \to 2$ of length~$9$, and gives
  $M_{\texttt{flow}}((1,2))=6$, against
  $M_{\texttt{top}}((1,2))=12$ from Corollary~\ref{cor:bigM-Lbar}. With $b=2$, one further arc
  of cost~$1$ can be interdicted together with $(1,2)$, and two paths
  that share no arc of cost~$1$ are needed: the paths~$1\to5\to3\to2$ and $1\to5\to4\to3\to2$ of
  \Cref{FIG_bigM_flow} give $\mu((1,2))=13$ and
  $M_{\texttt{flow}}((1,2))=13-\ell_{(1,2)}=10$.
  The arc~$(2,1)$ has cost~$2$, so that for $b=2$ it is interdicted
  alone and the single path~$2\to5\to4\to1$ of \Cref{FIG_bigM_flow1}
  suffices: $M_{\texttt{flow}}((2,1))=4$, against
  $M_{\texttt{top}}((2,1))=12$ from the arc-length rule. For $b=2$, the
  values of $M_{\texttt{rflow}}$ sum to~$191$ over all arcs and roots,
  against~$89$ for the thresholds~$\eta(a,w)$ of
  Example~\ref{ex:running-exact}.
\end{example}

\begin{table}[htb]
  \centering
  \caption{Penalty values of the running example of \Cref{FIG_a}: the
    root-independent rules $M_{\texttt{top}}$ and $M_{\texttt{flow}}$ for $b=1$ and
    $b=2$, and the root-dependent rules $M_{\texttt{rtop}}$ and $M_{\texttt{rflow}}$
    for $b=2$ and every root~$w$. A dash means that the arc cannot be
    interdicted, so that the penalty value is~$0$. For $b=2$, summing
    over all arcs and roots, the values are $520$ for $M_{\texttt{top}}$,
    $265$ for $M_{\texttt{flow}}$, $317$ for $M_{\texttt{rtop}}$, and $191$ for
    $M_{\texttt{rflow}}$.}
  \label{tab:running-bigM}
  \small
  \setlength{\tabcolsep}{3.5pt}
  \begin{tabular}{@{}crr r rr r rrrrr r rrrrr@{}}
    \toprule
    & \multicolumn{2}{c}{$M_{\texttt{top}}$} & & \multicolumn{2}{c}{$M_{\texttt{flow}}$} & & \multicolumn{5}{c}{$M_{\texttt{rtop}}$, $b=2$} & & \multicolumn{5}{c}{$M_{\texttt{rflow}}$, $b=2$}\\
    \cmidrule{2-3}\cmidrule{5-6}\cmidrule{8-12}\cmidrule{14-18}
    Arc & $b=1$ & $b=2$ & & $b=1$ & $b=2$ & & $w=1$ & $2$ & $3$ & $4$ & $5$ & & $w=1$ & $2$ & $3$ & $4$ & $5$\\
    \midrule
    $(1,2)$ & 12 & 12 & &  6 & 10 & & 12 &  0 &  5 &  9 &  6 & & 10 &  0 &  5 &  9 &  6\\
    $(1,5)$ & -- &  8 & & -- &  1 & &  8 &  8 &  5 &  8 &  0 & &  1 &  1 &  1 &  1 &  0\\
    $(2,1)$ & -- & 12 & & -- &  4 & &  0 & 12 &  6 &  6 &  6 & &  0 &  4 &  4 &  4 &  4\\
    $(2,5)$ & 13 & 13 & &  5 &  5 & &  9 & 13 &  9 &  7 &  0 & &  5 &  5 &  5 &  5 &  0\\
    $(3,2)$ & -- & 10 & & -- &  5 & &  6 &  0 & 10 &  8 &  9 & &  5 &  0 &  5 &  5 &  5\\
    $(3,4)$ & 11 & 11 & &  5 & 10 & &  5 &  7 & 11 &  0 &  8 & &  5 &  7 & 10 &  0 &  8\\
    $(4,1)$ & 11 & 11 & &  7 &  7 & &  0 &  8 &  7 & 11 &  9 & &  0 &  7 &  7 &  7 &  7\\
    $(4,3)$ & -- &  8 & & -- &  5 & &  4 &  6 &  0 &  8 &  7 & &  4 &  5 &  0 &  5 &  5\\
    $(5,3)$ & 10 & 10 & &  4 &  4 & &  8 & 10 &  0 &  5 & 10 & &  4 &  4 &  0 &  4 &  4\\
    $(5,4)$ & -- &  9 & & -- &  2 & &  7 &  9 &  6 &  0 &  9 & &  2 &  2 &  2 &  0 &  2\\
    \bottomrule
  \end{tabular}
\end{table}

All four rules are computable in polynomial time. The rule
$M_{\texttt{rtop}}$ only requires one APSP computation of the
non-interdicted graph in addition to sorting the arc lengths, whereas
the flow rules require minimum-cost flow computations. Since
$k_a\leq h_a$ and $k_0\leq h_0$, every flow has at most $m+1$ units,
so that it is computed by successive shortest paths and decomposed
into paths in time polynomial in the input size, independently of the
numerical value of the budget.
\Cref{tab:running-bigM} reports their values on the running example.
A step-by-step visual illustration of these rules on the
running example is available at
\url{https://fabiofurini.github.io/big-m-guide/}.

\subsection{Destination-Aware Poly-Time Computable Rules}
\label{sec:destination-bigM}

The four rules of \Cref{sec:bigM-bounds} ignore the requested
destination sets. A smaller penalty can suffice when a shortcut
through an interdicted arc does not reduce the distance to any
requested destination. For an active root~$w$ and a requested destination~$t$, let
$U(w,t)$ be a robust upper bound on the residual distance, i.e.,

\begin{equation}
  \label{eq:destination-upper-bound}
  U(w,t)\geq d^I(w,t)
  \quad\text{for every feasible interdiction set }I.
\end{equation}

The bound~$U(w,t)$ need not be exact. It only has to be large enough to
cover the distance that may result from any feasible interdiction. We use
the two poly-time computable choices already available in the paper. For every
ordered pair~$(p,t)$ of distinct vertices, set
\[
  U_{\texttt{top}}(p,t)
  \define
  \sigma(p,t),
  \quad
  U_{\texttt{flow}}(p,t)
  \define
  \min\Set{\mu(p,t),\sigma(p,t)},
\]
and set~$U_{\texttt{top}}(t,t)=U_{\texttt{flow}}(t,t)=0$.
Under Assumption~\ref{ass:feasibility}, both are robust upper bounds
for every requested pair~$(w,t)$: a shortest residual $w$--$t$ path is
simple and neither enters~$w$ nor leaves~$t$, so that its length is at
most $\sigma(w,t)$, and $d^I(w,t)\le\mu(w,t)$ by
Proposition~\ref{lem:disjoint-paths}. For an interdictable arc~$a=(u,v)$ and an active root~$w$, define
\begin{align}
  M_{\texttt{dtop}}(a,w)
  &\define
  \max_{t\in\mathcal T_w}
  \max\Set{
    0,
    U_{\texttt{top}}(w,t)-d(w,u)-\ell_a-d(v,t)
  },
  \label{eq:M-dtop}\\
  M_{\texttt{dflow}}(a,w)
  &\define
  \max_{t\in\mathcal T_w}
  \max\Set{
    0,
    U_{\texttt{flow}}(w,t)-d(w,u)-\ell_a-d(v,t)
  }.
  \label{eq:M-dflow}
\end{align}
That is, $M_{\texttt{dtop}}(a,w)$ and $M_{\texttt{dflow}}(a,w)$ are the
largest amount, over the requested destinations~$t$ of~$w$, by which
$U_{\texttt{top}}(w,t)$ and $U_{\texttt{flow}}(w,t)$, respectively,
exceed the length~$d(w,u)+\ell_a+d(v,t)$ of a shortest $w$--$t$ walk
through~$a$ in~$\mathcal G$.
For arcs with $c_a>b$, we set
$M_{\texttt{dtop}}(a,w)=M_{\texttt{dflow}}(a,w)=0$.

To interpret this formula, consider a $w$--$t$ path that uses an interdicted arc
$a=(u,v)$. Its penalized length is at least
\[
  d(w,u)+\ell_a+M(a,w)+d(v,t),
\]
because all arc lengths and penalties are nonnegative. Therefore, if
$M(a,w)$ is at least the corresponding value~$M_{\texttt{dtop}}(a,w)$ or $M_{\texttt{dflow}}(a,w)$, every requested
path that uses the interdicted arc~$a$ is at least as long as the
corresponding bound~$U_{\texttt{top}}(w,t)$ or $U_{\texttt{flow}}(w,t)$,
and hence at least as long as the residual shortest $w$--$t$ path.

The following proposition formalizes this argument.
\begin{proposition}[Destination-aware poly-time computable penalties]
  \label{prop:destination-polynomial}
  Under Assumption~\ref{ass:feasibility}, the values~$M_{\texttt{dtop}}(a,w)$ and $M_{\texttt{dflow}}(a,w)$ defined in
  \eqref{eq:M-dtop}--\eqref{eq:M-dflow} are valid penalties for
  the MPSP.
\end{proposition}
\begin{proof}
  The proof is identical for the two rules. Write
  $(M_\star,U_\star)$ for either
  $(M_{\texttt{dtop}},U_{\texttt{top}})$ or
  $(M_{\texttt{dflow}},U_{\texttt{flow}})$.
  Fix a feasible interdiction set~$I$, an active root~$w$, a requested
  destination~$t\in\mathcal T_w$, and a $w$--$t$ path~$P$, with
  penalized length~$\ell^{M_\star}_{I,w}(P)$. If $P$ contains no
  interdicted arc, it is a residual path and its length is at least
  $d^I(w,t)$. Otherwise, let $a=(u,v)\in I$ be an interdicted arc
  on~$P$. The prefix of~$P$ up to $u$ has length at least $d(w,u)$ and
  the suffix from~$v$ to~$t$ has length at least $d(v,t)$. Since the
  penalties are nonnegative, the definition of $M_\star$ gives
  \[
    \ell^{M_\star}_{I,w}(P)
    \ge
    d(w,u)+\ell_a+M_\star(a,w)+d(v,t)
    \ge
    U_\star(w,t)
    \ge
    d^I(w,t).
  \]
  Thus, every penalized $w$--$t$ path is at least as long as a residual
  shortest path. Since every residual path is still available in the
  penalized graph, the two shortest-path distances are equal.
\end{proof}

Both rules are poly-time computable. The
maximum is taken only over the destinations in~$\mathcal T_w$. A vertex that is not requested from~$w$ imposes no condition on the penalty.
The gain comes from the term $d(v,t)$: an interdicted arc only has to
be made unattractive for the destinations that are actually requested
from~$w$.

The destination-aware values can be combined with the baseline
rule~\eqref{eq:M-rflow} by keeping, for each arc and root, the smaller
of the two penalties. The following proposition states that this
combination is valid and coincides with the baseline when all the
destinations are requested.

\begin{proposition}[Hybrid destination-aware penalties]
  \label{prop:destination-hybrid}
  Under Assumption~\ref{ass:feasibility}, the values
  \[
    \min\set{M_{\texttt{rflow}}(a,w),M_{\texttt{dflow}}(a,w)}
  \]
  are
  valid penalties for the MPSP. Moreover, if
  $\mathcal T_w=\mathcal V\setminus\set{w}$ for an active root~$w$,
  then $\min\set{M_{\texttt{rflow}}(a,w),M_{\texttt{dflow}}(a,w)}
  =M_{\texttt{rflow}}(a,w)$ for every interdictable arc~$a$.
\end{proposition}
\begin{proof}
  Let $M$ denote the hybrid values, fix a feasible interdiction
  set~$I$, an active root~$w$, a requested destination~$t\in\mathcal T_w$, and a $w$--$t$ path~$P$. If $P$ contains an
  interdicted arc~$a=(u,v)$ with
  $M_{\texttt{dflow}}(a,w)\leq M_{\texttt{rflow}}(a,w)$, the penalty of~$a$
  is $M_{\texttt{dflow}}(a,w)$, and the chain of inequalities in the
  proof of Proposition~\ref{prop:destination-polynomial} gives
  $\ell^M_{I,w}(P)\geq d^I(w,t)$, irrespective of the other penalties,
  which are nonnegative. Otherwise, every interdicted arc of~$P$
  receives its penalty~$M_{\texttt{rflow}}(a,w)$, so that
  $\ell^M_{I,w}(P)=\ell^{M_{\texttt{rflow}}}_{I,w}(P)$, which is at least
  $d^I(w,t)$ by Proposition~\ref{prop:polynomial-rules}. Since every
  residual path is still available in the penalized graph, the two
  shortest-path distances are equal.

  For the second statement, let $a=(u,v)$ be interdictable. If $w=v$,
  then $M_{\texttt{rflow}}(a,w)=0\leq M_{\texttt{dflow}}(a,w)$. Otherwise, $v\in\mathcal T_w$ and,
  since $d(v,v)=0$, the term $t=v$ in~\eqref{eq:M-dflow} gives
  $M_{\texttt{dflow}}(a,w)\geq
  \max\set{0,\min\set{\mu(w,v),\sigma(w,v)}-d(w,u)-\ell_a}
  \geq M_{\texttt{rflow}}(a,w)$.
\end{proof}

\section{Supervalid Inequalities for the Compact Formulation}
\label{sec:svi}

In this section, we adapt the supervalid inequalities (SVIs)
of~\citet{Israeli_Wood:2002} to the MPSP interdiction problem and add
them to the compact Formulation~\eqref{eq:single-level} as lazy
constraints. No decomposition is involved and the inequalities are
generated from the interdiction decisions found during the
branch-and-bound search of Formulation~\eqref{eq:single-level}. The
basic inequality does not depend on the penalties, while its
cardinality version is valid only with valid penalties and becomes
stronger as the penalties decrease.

For a feasible interdiction decision~${\boldsymbol x}$, let
$I({\boldsymbol x})\define\set{a\in\mathcal A:x_a=1}$ and
\[
  F({\boldsymbol x})
  \define
  \varphi_{\mathcal K}\left(
    \mathcal G(\mathcal V,\mathcal A\setminus I({\boldsymbol x})),
    {\boldsymbol \ell}
  \right)
  =
  \sum_{(w,t)\in\mathcal K}d^{I({\boldsymbol x})}(w,t),
\]
i.e., $F({\boldsymbol x})$ is the MPSP value of the interdicted network.
Let $\bar f$ be the value of the best interdiction decision found so
far, called the \emph{incumbent value}. We call an inequality
\emph{supervalid} with respect to~$\bar f$ if it is satisfied by every
feasible interdiction decision~${\boldsymbol x}$ with
$F({\boldsymbol x})>\bar f$. This is stronger than Definition~1
of~\citet{Israeli_Wood:2002}, which only requires that the inequality
does not eliminate all optimal solutions unless the incumbent is
already optimal, and it is the property that their proofs establish.
Unlike the weaker one, it is preserved when several inequalities are
added together: any family of supervalid inequalities can be added to
Formulation~\eqref{eq:single-level} without removing an improving
decision, and, if no feasible decision satisfies them, $\bar f$ is the
optimal value.

Let $\hat{\boldsymbol x}$ be a feasible interdiction decision with
$F(\hat{\boldsymbol x})\le\bar f$, and let
$\hat I\define I(\hat{\boldsymbol x})$. For each active root~$w\in\mathcal R$, select one shortest $w$--$t$ path in the interdicted network $\mathcal G(\mathcal V,\mathcal A\setminus\hat I)$ for every
destination~$t\in\mathcal T_w$, and let $\hat y_{aw}$ be the number of
selected paths rooted in~$w$ that use arc~$a$. The vector~$\hat{\boldsymbol y}$ is feasible for the follower of the penalty
formulation~\eqref{eq:apspip-reform}, it uses no interdicted arc, and
\begin{equation}
  \label{eq:svi-value}
  \sum_{w\in\mathcal R}\sum_{a\in\mathcal A}\ell_a\hat y_{aw}
  =F(\hat{\boldsymbol x}).
\end{equation}
The \emph{support} of $\hat{\boldsymbol y}$ is the set of arcs used by
at least one selected path:
\begin{equation}
  \label{eq:svi-support}
  \operatorname{supp}(\hat{\boldsymbol y})
  \define
  \Set{
    a\in\mathcal A:
    \sum_{w\in\mathcal R}\hat y_{aw}>0
  }.
\end{equation}

The following proposition extends the basic SVI
of~\citet{Israeli_Wood:2002} to several pairs.

\begin{proposition}
  \label{prop:svi-basic}
  The inequality
  \begin{equation}
    \label{eq:svi-basic}
    \sum_{a\in\operatorname{supp}(\hat{\boldsymbol y})}x_a\;\ge\;1
  \end{equation}
  is supervalid with respect to~$\bar f$.
\end{proposition}
\begin{proof}
  Let ${\boldsymbol x}$ be a feasible interdiction decision that
  violates~\eqref{eq:svi-basic}. Since $I({\boldsymbol x})$ contains no
  arc of $\operatorname{supp}(\hat{\boldsymbol y})$, every selected path
  remains available in the interdicted network of~${\boldsymbol x}$. Hence,
  for every $(w,t)\in\mathcal K$, $d^{I({\boldsymbol x})}(w,t)$ is at
  most the length of the selected $w$--$t$ path. Summing over the
  requested pairs and using~\eqref{eq:svi-value}, we obtain
  $F({\boldsymbol x})\le F(\hat{\boldsymbol x})\le\bar f$.
\end{proof}

Inequality~\eqref{eq:svi-basic} states that a decision improving
on~$\bar f$ must interdict at least one arc of the selected paths.
The support is taken over all active roots: an improving decision need
not increase the contribution of every root. With a proper pair set,
the paths to vertices outside~$\mathcal T_w$ are not selected, so the
support can be smaller than in the all-pairs case, and the inequality
stronger. For a single pair~$(s,t)$, the support is one shortest $s$--$t$
path of the interdicted network, and Inequality~\eqref{eq:svi-basic} is the
covering constraint of the decomposition of~\citet{Israeli_Wood:2002}.

The penalties give a stronger inequality. For each arc~$a\in\mathcal A$, let
\[
  \hat q_a
  \define
  \sum_{w\in\mathcal R}M(a,w)\,\hat y_{aw},
\]
i.e., the total penalty that the selected paths pay on arc~$a$ if
$a$ is interdicted, and, for each integer $\rho$ with
$0\le\rho\le\card{\operatorname{supp}(\hat{\boldsymbol y})}$, let
$Q_\rho$ be the sum of the $\rho$ largest values~$\hat q_a$ over
$a\in\operatorname{supp}(\hat{\boldsymbol y})$, with $Q_0=0$. Thus,
$Q_\rho$ is an upper bound on the total penalty that the selected paths
pay when at most $\rho$ arcs of the support are interdicted.

The following proposition raises the right-hand side of Inequality~\eqref{eq:svi-basic} from one to~$\rho+1$.
\begin{proposition}
  \label{prop:svi-cardinality}
  Let the penalties~$M(a,w)$ be valid, and let $\rho$ be an integer
  with $0\le\rho\le\card{\operatorname{supp}(\hat{\boldsymbol y})}$
  such that $F(\hat{\boldsymbol x})+Q_\rho\le\bar f$. Then, the inequality
  \begin{equation}
    \label{eq:svi-cardinality}
    \sum_{a\in\operatorname{supp}(\hat{\boldsymbol y})}x_a\;\ge\;\rho+1
  \end{equation}
  is supervalid with respect to~$\bar f$.
\end{proposition}
\begin{proof}
  Let ${\boldsymbol x}$ be a feasible interdiction decision. Since the
  penalties are valid, $F({\boldsymbol x})$ is the optimal value of the
  follower of~\eqref{eq:apspip-reform} for~${\boldsymbol x}$. Since
  $\hat{\boldsymbol y}$ is feasible for this follower, we have
  \[
    F({\boldsymbol x})
    \le
    \sum_{w\in\mathcal R}\sum_{a\in\mathcal A}
    (\ell_a+M(a,w)x_a)\hat y_{aw}
    \overset{\eqref{eq:svi-value}}{=}
    F(\hat{\boldsymbol x})
    +\sum_{a\in\operatorname{supp}(\hat{\boldsymbol y})}\hat q_a x_a,
  \]
  where the sum can be restricted to the support because $\hat q_a=0$
  for every other arc. Assume that ${\boldsymbol x}$
  violates~\eqref{eq:svi-cardinality}, i.e., at most $\rho$ arcs of the
  support are interdicted. Since ${\boldsymbol x}$ is binary and
  $\hat q_a\ge0$, the last sum is at most~$Q_\rho$. Hence,
  $F({\boldsymbol x})\le F(\hat{\boldsymbol x})+Q_\rho\le\bar f$.
\end{proof}

For $\rho=0$, Inequality~\eqref{eq:svi-cardinality}
coincides with~\eqref{eq:svi-basic}. We always use the largest~$\rho$
that satisfies the condition of Proposition~\ref{prop:svi-cardinality}.
If the condition holds for
$\rho=\card{\operatorname{supp}(\hat{\boldsymbol y})}$, no feasible
decision satisfies~\eqref{eq:svi-cardinality}, which certifies
that~$\bar f$ is optimal, and the search stops. The cardinality inequality is the counterpart
of the tightening of~\citet[Corollary~1]{Israeli_Wood:2002}, with the
delays replaced by the penalties, and it links the inequalities to
\Cref{sec:bigM-analysis}. On the one hand, the proof requires valid
penalties: with penalties that are too small, the upper bound on
$F({\boldsymbol x})$ can fail and the inequality can remove improving
decisions. On the other hand, for fixed $\hat{\boldsymbol x}$,
$\hat{\boldsymbol y}$, and~$\bar f$, the coefficients $\hat q_a$ are
nondecreasing in the penalties, so that smaller valid penalties give a
value of~$\rho$ at least as large.

\subsection{Algorithmic Setup}

The inequalities are separated at every integer solution found during
the branch-and-bound search of Formulation~\eqref{eq:single-level}.
Given such a solution~$\hat{\boldsymbol x}$, we compute one shortest-path
tree of the interdicted network for each active root, with
$\card{\mathcal R}$ runs of Dijkstra's algorithm, which give
$\hat{\boldsymbol y}$ and $F(\hat{\boldsymbol x})$. We then set~$\bar f\define\max\set{\bar f,F(\hat{\boldsymbol x})}$ and add
Inequality~\eqref{eq:svi-cardinality} as a lazy constraint. Since
$\hat{\boldsymbol x}$ interdicts no arc of
$\operatorname{supp}(\hat{\boldsymbol y})$, every inequality cuts off
the solution that generated it. Hence, the solver never accepts an
incumbent of its own, and the best evaluated decision is optimal as
soon as no integer solution is left, or as soon as the upper bound of
the remaining search is at most~$\bar f$, since every excluded decision
has value at most~$\bar f$ by supervalidity. The search stops at the
first of these events. Since the solver accepts no incumbent, it cannot
prune the nodes whose bound does not exceed~$\bar f$. When the lengths
are integer, we therefore also consider a variant with an
\emph{objective cutoff}: every time $\bar f$ increases, we add the lazy
constraint stating that the objective function of
Formulation~\eqref{eq:single-level} is at least $\bar f+ 1/2$. This
constraint only excludes decisions of value at most~$\bar f$, i.e.,
it is the usual pruning by the incumbent value. In contrast,
\citet[Section~7.3]{Israeli_Wood:2002} generate the SVIs with a
heuristic before the solution of the compact formulation and add them
to it all at once, on instances with finite delays.

The basic inequalities alone also give an exact method without
penalties, which adapts the covering decomposition
of~\citet[Algorithm~2]{Israeli_Wood:2002} to the sum over the requested
pairs. Its master problem contains only the interdiction variables, the
budget constraint, and the basic inequalities~\eqref{eq:svi-basic}
generated so far. In~\citet{Israeli_Wood:2002}, the master problem is a
feasibility problem, and its objective function plays no role in the
correctness of the method: in our implementation it maximizes the
total interdiction cost~$\sum_{a\in\mathcal A}c_ax_a$, which steers the
search towards decisions that use the budget, and its value is not a
bound on the MPSP interdiction value. At every iteration,
a feasible solution of the master problem is evaluated by one
shortest-path tree for each active root, the incumbent is updated, and
the corresponding basic inequality is added. When the master problem
becomes infeasible, the incumbent is optimal by
Proposition~\ref{prop:svi-basic}. We solve the master problem in a
single branch-and-bound search, in which every integer solution is
evaluated and cut off by its basic inequality as a lazy constraint,
instead of solving a new master problem at every iteration, and we do
not use the greedy heuristic of~\citet{Israeli_Wood:2002} for the
master problem. With a single pair, the support of every inequality is
one shortest path, and the method generates the same inequalities as
Algorithm~2 of~\citet{Israeli_Wood:2002} with arc destruction. With several pairs, it is the
union of the selected paths of all requested pairs, which tends to
contain most arcs as the number of pairs grows. Each inequality then
excludes few interdictions, and the decomposition needs many
iterations.

\section{Computational Experiments}
\label{sec:computational}

In this section, we evaluate the effect of the penalties on the compact
Formulation~\eqref{eq:single-level}. The study pursues five goals: (i)
measuring how the penalty rules affect the linear optimization
relaxation, the branch-and-bound search, and the number of instances
solved to proven optimality, (ii) measuring how far the poly-time computable
rules are from the componentwise minimal penalties~$\eta(a,w)$, (iii) assessing whether supervalid inequalities, or the penalty-free
covering decomposition of \citet{Israeli_Wood:2002}, can replace
stronger penalties, (iv) assessing whether destination-aware
penalties pay off when only a subset of the pairs is requested, and
how the covering decomposition behaves as the number of pairs grows, and (v) comparing the methods in the single-pair setting
of~\citet{Israeli_Wood:2002}. Goals (i)--(iii) concern the all-pairs
specialization $\mathcal K=\mathcal K^{\mathrm{AP}}$, goal (iv) pair
sets~$\mathcal K\subsetneq\mathcal K^{\mathrm{AP}}$, and goal (v) a single pair.

\subsection{Implementation and Setup}
\label{sec:setup}

Formulation~\eqref{eq:single-level} is implemented in C++17 on top of
the C++ interface of \textit{Gurobi~13.0}~\citep{Gurobi:2025} and
compiled at optimization level \texttt{-O3}. The penalties are computed
before the model is built. Shortest paths are computed by our own
implementation of Dijkstra's algorithm with a binary heap, and the
minimum-cost flows of the surviving-path bounds by the successive
shortest path algorithm with node potentials. Two implementations of the formulation are
available: one specialized to $\mathcal K=\mathcal K^{\mathrm{AP}}$
and one for arbitrary pair sets. As a correctness check, the two
implementations, with $\mathcal K=\mathcal K^{\mathrm{AP}}$ in the
second one, and the enumeration of the feasible interdiction sets
return the same optimal value on all $108$ small instances.

Shortest-path computations are needed in two places only: before
the model is built, to compute the distances used by the penalty rules,
and, in the configurations with supervalid inequalities and in the
covering decomposition, to evaluate every candidate interdiction.
Formulation~\eqref{eq:single-level} itself requires none during the
branch-and-bound search. In both cases, the requested pairs are grouped
by origin, so that a single shortest-path tree from each active root
yields the distances to all its requested destinations, and an
evaluation requires at most $\card{\mathcal R}$ runs of Dijkstra's
algorithm. On our small and scaling instances, which have at most $100$ vertices and
$1412$ arcs, these computations have a negligible impact on the
running times: for instance, the complete computation of the penalties
$M_{\texttt{rflow}}$, minimum-cost flows included, has a median computation
time of $0.17$~seconds and a maximum of $4.2$~seconds. A standard implementation of
Dijkstra's algorithm is therefore sufficient. On large graphs, in
particular when few destinations per origin are requested, this effort
can be reduced by bidirectional search or by dedicated multiple-pairs
shortest-path algorithms such as those of~\citet{Grappe_etal}.

{
All runs were executed on a virtual machine with an Intel Xeon Gold 5418Y (10 cores)
and 58~GB of RAM under Ubuntu 22.04.
}
Gurobi runs in single-threaded mode with a relative optimality gap of~$0$ and all
other parameters at their default values. Six runs are executed in
parallel, at low priority, and no memory limit is imposed. Since the runs share the machine, the reported times are larger than
the times of sequential runs, and differences of a few percent between
configurations are within the variability of this setting.
{
A time limit of $1800$~seconds is imposed on the solver for each
run. The time needed to compute the penalties is reported separately
and is not counted in the time limit. Unless stated otherwise, the
\emph{time} of a run is the time to build the model and to solve it,
without the computation of the penalties. In \Cref{tab:e5,tab:iw},
which compare the compact formulation with the covering decomposition,
the time of the compact formulation also includes the computation of
the penalties, and so do the times of \Cref{fig:e1-solved}. In \Cref{tab:e3}, both timings are reported, with total
method time including penalty computation.}

{The code, the instances, and the results of all runs
are available at \url{https://doi.org/10.5281/zenodo.24680030}.}

For each instance, let $z^{\mathrm{best}}$ be the best known value,
i.e., the optimal value when at least one run proves it, and the best
value found by any run otherwise. For a run, let $z^{\mathrm{LP}}$ be
the optimal value of the linear optimization relaxation of
Formulation~\eqref{eq:single-level}, $z^{\mathrm{root}}$ the upper
bound at the end of the root node, and $\bar z$ and $\underline z$ the
best upper bound and the value of the best interdiction found within
the time limit. The \emph{LP gap} and the \emph{root gap} are
$100\,(z^{\mathrm{LP}}-z^{\mathrm{best}})/z^{\mathrm{best}}$ and
$100\,(z^{\mathrm{root}}-z^{\mathrm{best}})/z^{\mathrm{best}}$, in
percent, and measure the strength of the formulation independently of
the search. When no run proves optimality, they are measured with
respect to the best known value rather than to the optimum. The root
bound includes the presolve and the cutting planes of the solver. The \emph{final gap} of a run that reaches the time limit
is $100\,(\bar z-\underline z)/\underline z$.
Search times and node counts are geometric means over
the instances solved by the corresponding approach.
For penalty-computation times in \Cref{tab:e1}, we use a geometric
mean over all instances, with each time floored at
$10^{-4}$~seconds. A run that
reaches the time limit is excluded from the search-time and node-count
means, but contributes to the LP, root, and final-gap statistics.
The LP
and root gaps are arithmetic means, and the final gaps are medians over
the runs that reach the time limit.
Runtime means are conditional on each approach's solved set. When
approaches solve different numbers of instances, the means should therefore
be interpreted jointly with the solved counts rather than as a standalone
ranking.

\subsection{Benchmark Instances}
\label{sec:instances}

We use two sets of instances, generated by the same procedure. All of
them satisfy Assumption~\ref{ass:strong-residual}, i.e., no
interdiction within the budget disconnects the graph.

The \emph{small} set contains $108$ instances on which the minimal
penalties~$\eta(a,w)$ and the optimal values can be computed by
enumerating the feasible interdiction sets. It consists of
Erd\H{o}s--R\'enyi graphs with $n_0\in\set{6,8,10,12}$ vertices before
the extraction described below, edge probability $p\in\set{0.5,0.7}$,
unit costs, lengths in $\set{1}$ or $\set{1,\dots,10}$, and
$b\in\set{1,2,3}$, Erd\H{o}s--R\'enyi graphs with $n_0\in\set{8,10}$,
$p=0.8$, lengths in $\set{1,\dots,10}$, costs in $\set{1,2,3}$, and
$b\in\set{2,3,4}$, toroidal grids of sizes $3\times3$, $3\times4$, and
$4\times4$ with unit costs and $b\in\set{1,2,3}$, and three hand-made
examples.

The \emph{scaling} set contains $229$ instances with
$n_0\in\set{20,30,40,50,60,80,100}$:
Erd\H{o}s--R\'enyi graphs with average degree~$10$ (ER) and
Barab\'asi--Albert graphs with~$6$ edges per new vertex (BA), with unit
costs and lengths in $\set{1}$ or $\set{1,\dots,10}$, Erd\H{o}s--R\'enyi graphs with average degree~$14$, lengths in
$\set{1,\dots,10}$, and costs in $\set{1,2,3}$ (ER, costs), random geometric graphs with $n_0$ points in the unit square and radius
$\sqrt{13/(\pi n_0)}$, i.e., expected average degree about~$13$, in
which the length of an edge is $10$~times its Euclidean length divided
by the radius, rounded up, and lies in $\set{1,\dots,10}$ (only
$n_0\ge30$ for these graphs), and toroidal grids of sizes $5\times5$ to $10\times10$ without (torus)
and with diagonals (king torus). For each random family, one seed is
used. A second seed gives $140$ further instances of the random
families, used only to check the robustness of the results on the
penalty rules. Of each graph we keep the largest $(b_{\max}+1)$-edge-connected
subgraph, where $b_{\max}$ is the largest budget of its family
(\texttt{k\_edge\_subgraphs} of NetworkX~\citep{Hagberg_etal:2008}),
and we replace every edge by two anti-parallel arcs, whose lengths are
drawn independently except for the geometric graphs. For every graph
and budget, Assumption~\ref{ass:strong-residual} is certified by
computing the minimum directed cut with the interdiction costs as
capacities. By Proposition~\ref{obs:cut}, the instance is kept only if
this cut exceeds~$b$. \Cref{tab:instances} summarizes the two sets.

\begin{table}[htb]
  \centering
  \caption{Benchmark instances: number of instances and ranges of the
    numbers of vertices~$n$ and arcs~$m$ (after the extraction of the
    edge-connected subgraph) and of the budget~$b$.}
  \label{tab:instances}
  \small
  \begin{tabular}{@{}llrccc@{}}
    \toprule
    Set & Family & Instances & $n$ & $m$ & $b$\\
    \midrule
    small   & ER         &  78 & 5--12   & 20--92    & 1--3\\
            & ER, costs  &   9 & 8--10   & 46--74    & 2--4\\
            & torus      &  18 & 9--16   & 36--64    & 1--3\\
            & hand-made  &   3 & 3--4    & 6--10     & 0--1\\
    \cmidrule{2-6}
            & Total      & 108 & 3--16   & 6--92     & 0--4\\
    \midrule
    scaling & ER         &  56 & 20--92  & 206--954  & 2--5\\
            & BA         &  56 & 16--100 & 134--1128 & 2--5\\
            & ER, costs  &  21 & 20--100 & 282--1412 & 2, 4, 6\\
            & geometric  &  24 & 22--85  & 190--988  & 2--5\\
            & torus      &  24 & 25--100 & 100--400  & 2, 3\\
            & king torus &  48 & 25--100 & 200--800  & 2--5\\
    \cmidrule{2-6}
            & Total      & 229 & 16--100 & 100--1412 & 2--6\\
    \bottomrule
  \end{tabular}
\end{table}

\subsection{Penalty Rules and Methods}
\label{sec:rules-compared}

We compare five penalty rules. The first is the uniform penalty
\[
  M_{\texttt{unif}}\define n \, \max_{a\in\mathcal A}\ell_a,
\]
which is the sufficiently large delay proposed
by~\citet[Section~2]{Israeli_Wood:2002} to model arc destruction. It is valid, since every simple path has at most $n-1$
arcs. The other four are the poly-time computable rules $M_{\texttt{top}}$,
$M_{\texttt{flow}}$, $M_{\texttt{rtop}}$, and $M_{\texttt{rflow}}$ of
\Cref{sec:bigM-summary}, each defined as the minimum of the bounds
available at its level, so that
$M_{\texttt{rflow}}\le M_{\texttt{flow}}\le M_{\texttt{top}}$ and
$M_{\texttt{rflow}}\le M_{\texttt{rtop}}\le M_{\texttt{top}}$
componentwise. For every rule, the penalty of an arc with $c_a>b$ is
set to~$0$.
The rules differ in their computing cost: $M_{\texttt{unif}}$ is
computed in $O(m)$ time, $M_{\texttt{top}}$ requires sorting the arc
lengths, $M_{\texttt{rtop}}$ additionally requires one APSP
computation, $M_{\texttt{flow}}$ one minimum-cost flow per
interdictable arc, and $M_{\texttt{rflow}}$, in addition, one
minimum-cost flow per ordered pair of vertices. For pair sets, we also
use the \emph{destination-aware hybrid penalties}
$\min\set{M_{\texttt{rflow}},M_{\texttt{dflow}}}$, with
$M_{\texttt{dflow}}$ the rule of~\eqref{eq:M-dflow}. They are valid by
Proposition~\ref{prop:destination-hybrid}.

The penalty rules are combined with four algorithmic configurations,
which the tables name as follows:
\begin{itemize}
  \item \emph{compact}: Formulation~\eqref{eq:single-level} with one of
    the penalty rules, solved by the MIP solver. In a table or figure in
    which every row or curve uses it, a row or curve is named only by
    its penalty rule.
  \item \emph{compact + cardinality SVIs}: the compact formulation with
    the cardinality inequalities~\eqref{eq:svi-cardinality} added as
    lazy constraints, as described in \Cref{sec:svi}.
  \item \emph{compact + cardinality SVIs + cutoff}: the previous
    configuration with the objective cutoff of \Cref{sec:svi}.
  \item \emph{covering (basic SVIs)}: our adaptation of the covering
    decomposition of~\citet{Israeli_Wood:2002} to the sum over the
    requested pairs, described in \Cref{sec:svi}. Its master problem
    contains only the basic inequalities~\eqref{eq:svi-basic}, which are
    its covering constraints, and it uses no penalties.
\end{itemize}

\subsection{Effect of the Penalties on the Compact Formulation}
\label{sec:results-penalties}

\Cref{tab:e1} reports, for each rule, the results of
Formulation~\eqref{eq:single-level} on the $229$ scaling instances:
  the number of instances solved to proven optimality, the
  geometric means of the search time and of the number of nodes over the
  instances solved by each rule, the mean LP gap
and root gap, the median final gap of the runs that reach the time
limit, the mean penalty relative to $M_{\texttt{unif}}$, the share of
zero penalties, and the geometric mean of the time spent computing the
penalties.

\begin{table}[htb]
  \centering
  \caption{Compact formulation with the five penalty rules on the $229$
    scaling instances: instances solved to proven optimality,
    geometric mean of the time (seconds) and of the number of nodes
    over the instances solved by the corresponding rule,
    arithmetic mean of the LP gap and of the root gap (\%), median final
    gap (\%) of the runs that reach the time limit, mean penalty as a
    percentage of $M_{\texttt{unif}}$, share (\%) of zero penalties, and
    geometric mean of the time (seconds) needed to compute the
    penalties, which is not counted in the time limit. The relative
    mean penalty and the share of zero penalties are computed for each
    instance over the entries $(a,w)$ with $c_a\le b$ and $w\neq v$,
    where $a=(u,v)$, and then averaged arithmetically with equal weight
    for each instance. The LP gap
    is averaged over the $196$ instances on which the relaxation is
    solved within $60$~seconds for every rule, and the root gap over the
    $224$ instances with a root bound for every rule. Time limit:
    $1{,}800$~seconds. Numbers in parentheses after solved counts are
    percentages of the total instances in the corresponding block.}
  \label{tab:e1}
  \small
  \setlength{\tabcolsep}{4pt}
  {
\begin{tabular}{@{}l rrr r rrr r rrr@{}}
  \toprule
  & \multicolumn{3}{c}{search} & & \multicolumn{3}{c}{gap (\%)} & &
    \multicolumn{3}{c}{penalties}\\
  \cmidrule{2-4}\cmidrule{6-8}\cmidrule{10-12}
  Rule & solved & time & nodes & & LP & root & final & &
    mean & zero & time\\
  \midrule
  $M_{\texttt{unif}}$ & 52 (22.7\%) & 99.4 & 1282 &  & 102.5 & 78.8 & 48.9 &  & 100.0 & 0.0 & 0.000\\
  $M_{\texttt{top}}$ & 49 (21.4\%) & 83.5 & 895 &  & 98.2 & 74.7 & 46.1 &  & 94.9 & 0.0 & 0.003\\
  $M_{\texttt{rtop}}$ & 52 (22.7\%) & 90.2 & 897 &  & 97.4 & 74.1 & 46.4 &  & 91.9 & 0.0 & 0.014\\
  $M_{\texttt{flow}}$ & 171 (74.7\%) & 15.8 & 37 &  & 5.9 & 1.0 & 0.6 &  & 2.7 & 7.2 & 0.032\\
  $M_{\texttt{rflow}}$ & 190 (83.0\%) & 8.6 & 24 &  & 4.6 & 0.7 & 0.5 &  & 0.8 & 55.7 & 0.181\\
  \bottomrule
\end{tabular}

}
\end{table}

{
The penalties determine which instances can be solved. Formulation~\eqref{eq:single-level}
with $M_{\texttt{rflow}}$ solves $190$ of the $229$ instances, against
$171$ with $M_{\texttt{flow}}$, $52$ with $M_{\texttt{unif}}$ and
$M_{\texttt{rtop}}$, and $49$ with $M_{\texttt{top}}$.
The geometric mean of the search time, computed only on
instances solved by the corresponding rule, is $99.4$ seconds with
$M_{\texttt{unif}}$ and $8.6$ seconds with $M_{\texttt{rflow}}$.
These conditional means must be read together with the solved counts:
$M_{\texttt{unif}}$ and $M_{\texttt{rflow}}$ solve $52$ and $190$ of the
$229$ instances, respectively. On
the $196$ instances on which every relaxation is solved, the mean LP
gap decreases from $102.5\%$ with $M_{\texttt{unif}}$ to $4.6\%$ with
$M_{\texttt{rflow}}$; on the $224$ common root-bound instances, the
mean root gap decreases from $78.8\%$ to $0.7\%$. The advantage of the surviving-path rules
persists on the $140$ instances of the second seed:
$M_{\texttt{rflow}}$ solves $121$, $M_{\texttt{flow}}$ solves $108$,
and the other three rules solve $25$--$27$ instances.

Of the two refinements of $M_{\texttt{top}}$, the surviving paths
account for most of the gain. Root dependence alone changes little:
$M_{\texttt{rtop}}$ keeps $91.9\%$ of the mean penalty of
$M_{\texttt{unif}}$. Bounding residual distances by surviving paths
reduces the mean penalty to $2.7\%$ of $M_{\texttt{unif}}$ and the LP
gap to $5.9\%$. On top of the surviving paths, root dependence pays
off: $M_{\texttt{rflow}}$ sets $55.7\%$ of the entries to~$0$, against
$7.2\%$ for $M_{\texttt{flow}}$, and solves $19$ additional instances;
$M_{\texttt{flow}}$ solves none that $M_{\texttt{rflow}}$ does not.
The instances that remain unsolved with $M_{\texttt{rflow}}$ are mostly
the largest ones: it solves $87$ of the $90$ instances with $n<40$,
against $35$ of the $55$ instances with $n\ge80$. Their median final
gap is $0.5\%$.

Penalty-computation times are small in absolute terms:
$M_{\texttt{rflow}}$ has a median computation time of $0.17$~seconds and a
maximum of $4.2$~seconds. Overall, the strength of the penalties decides whether
the compact formulation is a practical method. This agrees with the
effect that \citet{Lim_Smith:2007} observe when they refine the
penalties of a flow interdiction model, here for an objective that sums
shortest-path distances.}

{
\Cref{fig:e1-solved} plots, for each rule, the number of instances
solved within time~$t$ on the $229$ scaling instances. The $x$-axis
reports the time in seconds, on a logarithmic scale, and the $y$-axis
the percentage of instances solved to proven optimality.}

\begin{figure}[htb]
  \centering
  \input{fig_e1_solved.tex}
  \caption{Percentage of the $229$ scaling instances solved to
    proven optimality within time~$t$ (logarithmic scale) by the compact
    formulation, for each of the five penalty rules. Time limit:
    $1{,}800$~seconds. The times include the computation of the
    penalties, which takes at most $4.2$~seconds (a median of
    $0.17$~seconds for $M_{\texttt{rflow}}$). The curves of $M_{\texttt{unif}}$, $M_{\texttt{top}}$, and
    $M_{\texttt{rtop}}$ almost coincide.}
  \label{fig:e1-solved}
\end{figure}

{
On the scaling instances, the curves of $M_{\texttt{unif}}$,
$M_{\texttt{top}}$, and $M_{\texttt{rtop}}$ remain far below those of
the two surviving-path rules, with $M_{\texttt{rflow}}$ above
$M_{\texttt{flow}}$ over the whole time range. Within one second,
$M_{\texttt{rflow}}$ solves $36$ instances and $M_{\texttt{flow}}$
solves $21$; within ten seconds, the respective counts are $104$ and
$73$.}

\subsection{Distance from the Minimal Penalties}
\label{sec:results-eta}

On the small instances, the componentwise minimal penalties~$\eta(a,w)$ of Corollary~\ref{prop:rooted-bigM} are computed by
enumerating the feasible interdiction sets. They are used only as a
reference, to measure how much room is left between the poly-time computable
rules and the smallest valid penalties. For each rule, we consider all
entries $(a,w)$ with $c_a\le b$ and $w\neq v$, where $a=(u,v)$.
\Cref{tab:e2} reports the share of entries in which the rule attains
$\eta(a,w)$ and, over the entries with $\eta(a,w)>0$, statistics of
the ratio $M(a,w)/\eta(a,w)$.

\begin{table}[htb]
  \centering
  \caption{Penalty rules versus the minimal penalties~$\eta(a,w)$ on
    the $108$ small instances: share (\%) of entries $(a,w)$ with
    $M(a,w)=\eta(a,w)$, mean, median, and maximum of the ratio
    $M(a,w)/\eta(a,w)$ over the entries with $\eta(a,w)>0$, and ratio of the sums $\sum M(a,w)/\sum\eta(a,w)$ over all entries.}
  \label{tab:e2}
  \small
  
\begin{tabular}{@{}l r r rrr r@{}}
  \toprule
  & & & \multicolumn{3}{c}{$M(a,w)/\eta(a,w)$} & \\
  \cmidrule{4-6}
  Rule & $=\eta$ (\%) & & mean & median & max & $\sum M/\sum\eta$\\
  \midrule
  $M_{\texttt{unif}}$ & 0.0 &  & 31.98 & 16.67 & 2997.00 & 65.53\\
  $M_{\texttt{top}}$ & 0.0 &  & 25.38 & 13.50 & 136.00 & 52.17\\
  $M_{\texttt{rtop}}$ & 0.0 &  & 24.18 & 13.00 & 130.00 & 48.64\\
  $M_{\texttt{flow}}$ & 22.3 &  & 2.34 & 1.25 & 23.00 & 3.66\\
  $M_{\texttt{rflow}}$ & 86.6 &  & 1.26 & 1.00 & 12.00 & 1.36\\
  \bottomrule
\end{tabular}

\end{table}

The $108$ small instances have $53{,}341$ entries, and $\eta(a,w)=0$
on $74\%$ of them: for most arcs and roots, the penalty~$0$ is
valid. The rules that ignore the
surviving paths almost never detect this: $M_{\texttt{unif}}$,
$M_{\texttt{top}}$, and $M_{\texttt{rtop}}$ attain $\eta(a,w)$ on~$0$,
$7$, and~$10$ entries, respectively, and exceed it by a median factor
between~$13$ and~$17$ on the entries with $\eta(a,w)>0$. $M_{\texttt{flow}}$ attains $\eta(a,w)$ on
$22.3\%$ of the entries, and $M_{\texttt{rflow}}$ on $86.6\%$, with a
median ratio of~$1$ and a ratio of the sums of~$1.36$. Since
$\eta(a,w)=0$ on most entries, the agreement on the $13{,}771$ entries
with $\eta(a,w)>0$ is more informative: $M_{\texttt{rflow}}$ attains
$\eta(a,w)$ on $76.7\%$ of them, against $47.0\%$ for
$M_{\texttt{flow}}$. On the entries with $\eta(a,w)=0$,
$M_{\texttt{rflow}}$ is zero on $90.0\%$. Hence, on these
instances, $M_{\texttt{rflow}}$ is close to the smallest valid
penalties, and most of the room left by the other rules is closed by
combining root dependence with surviving paths.

The approximation quality also varies with the
budget and the range of the arc lengths. With unit lengths,
$M_{\texttt{rflow}}$ attains $\eta(a,w)$ on $99.8\%$ of the entries,
and with lengths in $\set{1,\dots,10}$ on $75.4\%$. For
$b=1,2,3$, the share is $97.4\%$, $89.2\%$, and $78.3\%$, and the ratio
of the sums is $1.11$, $1.30$, and $1.38$. For $b=4$, which occurs
only on the instances with costs in $\set{1,2,3}$, the share drops to
$41.3\%$ and the ratio of the sums is~$3.23$. Restricted to the
entries with $\eta(a,w)>0$, the share is $97.6\%$, $82.0\%$, $69.6\%$,
and $22.0\%$ for $b=1,2,3,4$, so that the approximation becomes
clearly less tight as the budget grows. A likely reason is
that, with nonunit costs, the cardinality cap $k_a$ of the surviving-path bounds is computed with
the smallest cost and is therefore weaker.

\subsection{Supervalid Inequalities and the Covering Decomposition}
\label{sec:results-svi}

For arc destruction, \citet{Israeli_Wood:2002} recommend a covering
decomposition whose master problem contains no penalties, strengthened
by supervalid inequalities (SVIs). As described in \Cref{sec:svi}, the
SVIs can also be added to Formulation~\eqref{eq:single-level}. In this section, we assess whether
SVIs or a penalty-free decomposition can replace stronger penalties.
We compare seven configurations: Formulation~\eqref{eq:single-level}
with $M_{\texttt{top}}$ and with $M_{\texttt{rflow}}$, without SVIs,
with SVIs, and with SVIs and the objective cutoff, and the covering
decomposition of~\citet{Israeli_Wood:2002}, adapted to the
sum over the pairs as described in \Cref{sec:svi}. The compact
formulation uses the cardinality inequalities~\eqref{eq:svi-cardinality},
and the covering decomposition the basic
inequalities~\eqref{eq:svi-basic}. The comparison uses the $108$ small
instances and the $124$ scaling instances with $n_0\le 50$. The time limit is the same for all configurations and, for the covering
decomposition, covers the whole algorithm.

\Cref{tab:e5} reports, for each configuration and instance set, the
    number of instances solved to proven optimality, the
    geometric mean of the time over the instances solved by the corresponding
    configuration, the median final gap of the runs that reach the time
    limit, and the geometric mean of the number of
    interdictions evaluated over the instances solved by the corresponding
    configuration.

\begin{widetable}[htb]
  \centering
  \caption{Compact formulation with $M_{\texttt{top}}$ and
    $M_{\texttt{rflow}}$, without SVIs, with SVIs, and with SVIs and
    the objective cutoff, and covering decomposition
    of~\citet{Israeli_Wood:2002} (basic SVIs), on the $108$ small instances and on
    the $124$ scaling instances with $n_0\le50$: instances solved to
    proven optimality, geometric mean of the time (seconds)
    over the instances solved by the corresponding configuration, median
    final gap (\%) of the runs that reach the time limit, and geometric
    mean of the number of evaluated interdictions over the
    instances solved by the corresponding configuration.
    The covering decomposition has no upper bound before it terminates
    ($\infty$). Time limit: $1{,}800$~seconds. Numbers in
    parentheses after solved counts are percentages of the total instances
    in the corresponding block.}
  \label{tab:e5}
  \small
  \setlength{\tabcolsep}{4pt}
  {
\begin{tabular}{@{}l rrrr r rrrr@{}}
  \toprule
  & \multicolumn{4}{c}{small ($108$)} & &
    \multicolumn{4}{c}{scaling, $n_0\le50$ ($124$)}\\
  \cmidrule{2-5}\cmidrule{7-10}
  Configuration & solved & time & gap & eval. & & solved & time & gap & eval.\\
  \midrule
  compact, $M_{\texttt{top}}$ & 108 (100.0\%) & 0.17 & -- & -- &  & 43 (34.7\%) & 88.20 & 32.6 & --\\
  compact, $M_{\texttt{top}}$ + cardinality SVIs & 105 (97.2\%) & 0.71 & 2.6 & 198 &  & 23 (18.5\%) & 84.06 & 32.3 & 832\\
  compact, $M_{\texttt{top}}$ + cardinality SVIs + cutoff & 108 (100.0\%) & 0.27 & -- & 8 &  & 55 (44.4\%) & 136.30 & 36.9 & 10\\
  compact, $M_{\texttt{rflow}}$ & 108 (100.0\%) & 0.02 & -- & -- &  & 116 (93.5\%) & 3.39 & 1.0 & --\\
  compact, $M_{\texttt{rflow}}$ + cardinality SVIs & 108 (100.0\%) & 0.04 & -- & 23 &  & 96 (77.4\%) & 17.05 & 0.5 & 287\\
  compact, $M_{\texttt{rflow}}$ + cardinality SVIs + cutoff & 108 (100.0\%) & 0.03 & -- & 7 &  & 113 (91.1\%) & 5.50 & 0.6 & 9\\
  covering (basic SVIs) & 99 (91.7\%) & 0.41 & $\infty$ & 434 &  & 13 (10.5\%) & 116.69 & $\infty$ & 5391\\
  \bottomrule
\end{tabular}

}
\end{widetable}

{
On the scaling subset, Formulation~\eqref{eq:single-level} with
$M_{\texttt{rflow}}$ solves $116$ of the $124$ instances, with a
  geometric mean time of $3.39$ seconds. With SVIs, it solves~$96$
instances, and~$113$ when the objective cutoff is added. The compact
formulation with $M_{\texttt{top}}$ solves between~$23$ and~$55$
instances in its three configurations, and the covering decomposition
solves~$13$. On the scaling subset, the solved counts therefore
range from $13$ for covering to $116$ for the compact formulation with
$M_{\texttt{rflow}}$ without SVIs; the conditional time and evaluation
means below should be read together with this difference. On the small set, the compact formulation with
$M_{\texttt{rflow}}$ solves all $108$ instances in every configuration;
the compact formulation with $M_{\texttt{top}}$ and SVIs solves~$105$,
and the covering decomposition solves~$99$.

Without the objective cutoff, our lazy separation of the SVIs
slows down the compact formulation with both rules. Since every SVI cuts off the integer
solution that generated it, the solver accepts no incumbent and cannot
prune the nodes whose bound does not exceed~$\bar f$. It evaluates and
cuts off many interdictions of lower value: on the scaling subset, the
  geometric mean is $287$ with $M_{\texttt{rflow}}$ and $832$
  with $M_{\texttt{top}}$, against~$9$ and~$10$, respectively, with the
cutoff. The cutoff improves both SVI configurations. With
$M_{\texttt{top}}$, adding SVIs and the cutoff increases the number of
solved scaling instances from $43$ to $55$. With $M_{\texttt{rflow}}$,
however, the configuration without SVIs remains better, solving $116$
  instances compared with $113$. Neither do the SVIs close the gap between
the two rules: with $M_{\texttt{top}}$, the compact formulation solves
at most~$55$ scaling instances, against~$116$ with
$M_{\texttt{rflow}}$ and no SVIs. The covering decomposition evaluates
  about $5{,}400$ interdictions on the scaling subset, in the geometric
  mean, since its master problem gives no bound on the MPSP
interdiction value before it terminates.}

Overall, on the tested all-pairs instances, neither the tested
lazy-SVI configurations nor the
penalty-free covering decomposition replace the strengthened
penalties. This is the opposite of the ranking that
\citet{Israeli_Wood:2002} observe for a single pair.

\subsection{Destination-Aware Penalties for Multiple Pairs}
\label{sec:results-mpsp}

To isolate the effect of the destination-aware rules of
\Cref{sec:destination-bigM}, we fix the number of active roots and
vary only the number of requested destinations per root. Since
Formulation~\eqref{eq:single-level} has one dual block per active root,
the model has the same variables and constraints for every pair set
with the same roots, and only the penalty coefficients and the
objective change. We use the scaling instances with one budget per
graph, $b=3$, or $b=4$ for the graphs with costs in $\set{1,2,3}$: the
$41$ instances with $n_0\le60$, with $\card{\mathcal R}=5$ active
roots, and the $20$ instances with $70\le n_0\le100$, with
$\card{\mathcal R}=10$. The active roots are drawn uniformly at random,
and, for each root, a random order of the other vertices is fixed. $\mathcal T_w$ consists of the first $1$, $\lceil0.1(n-1)\rceil$,
$\lceil0.25(n-1)\rceil$, $\lceil0.5(n-1)\rceil$, or $n-1$ vertices of
this order, so that the pair sets of an instance are nested. The baseline penalties
are $M_{\texttt{rflow}}$, and the destination-aware penalties are
$\min\set{M_{\texttt{rflow}},M_{\texttt{dflow}}}$, which are valid by
Proposition~\ref{prop:destination-hybrid}.

\Cref{tab:e3} reports, for each number of destinations per root, the
results with the baseline and with the destination-aware penalties,
together with the mean share of entries $(a,w)$ in which the
destination-aware penalty is strictly smaller than the baseline one.

\begin{widetable}[htb]
  \centering
  \caption{Compact formulation with the baseline penalties~$M_{\texttt{rflow}}$ and with the destination-aware hybrid penalties~$\min\set{M_{\texttt{rflow}},M_{\texttt{dflow}}}$ (hybrid), and
    covering decomposition of~\citet{Israeli_Wood:2002} adapted to
    several pairs (covering, basic SVIs), by number~$\card{\mathcal T_w}$ of
    requested destinations per root ($n-1$: all other vertices): mean share
    (\%) of entries with a strictly smaller destination-aware penalty,
    instances solved to proven optimality, geometric means
    of model-building/search time (time) and total method time (total),
    in seconds, computed over the instances solved by the corresponding
    approach, and mean LP gap (\%). Reduction shares are
    computed within each instance over entries $(a,w)$ with
    $a=(u,v)$, $c_a\le b$, and $w\in\mathcal R\setminus\Set{v}$,
    and then averaged arithmetically with equal weight for each
    instance. For the compact formulation,
    total includes penalty computation, while time excludes it. For
    covering, only total method time is shown. The baseline
    penalties are computed for all $n$ roots, the hybrid ones only for
    the active roots, so that the total times of the two compact
    configurations are not directly comparable. Time limit:
    $1{,}800$~seconds. Numbers in parentheses after solved counts are
    percentages of the total instances in the corresponding block.}
  \label{tab:e3}
  \small
  \setlength{\tabcolsep}{3pt}
  {
\begin{tabular}{@{}l r r rrrr r rrrr r rr@{}}
  \toprule
  & & & \multicolumn{4}{c}{compact, $M_{\texttt{rflow}}$} & &
    \multicolumn{4}{c}{compact, hybrid} & &
    \multicolumn{2}{c}{covering (basic SVIs)}\\
  \cmidrule{4-7}\cmidrule{9-12}\cmidrule{14-15}
  $\card{\mathcal T_w}$ & reduced (\%) & & solved & time & total & LP gap & &
    solved & time & total & LP gap & & solved & total\\
  \midrule
  \multicolumn{15}{@{}l}{$n_0\le60$, $41$ instances, $\card{\mathcal R}=5$}\\
  $1$ & 38.2 &  & 41 (100.0\%) & 0.12 & 0.201 & 13.8 &  & 41 (100.0\%) & 0.10 & 0.122 & 9.1 &  & 41 (100.0\%) & 0.12\\
  $10\%$ & 31.1 &  & 41 (100.0\%) & 0.14 & 0.215 & 9.9 &  & 41 (100.0\%) & 0.12 & 0.142 & 7.4 &  & 41 (100.0\%) & 4.27\\
  $25\%$ & 21.9 &  & 41 (100.0\%) & 0.17 & 0.255 & 8.6 &  & 41 (100.0\%) & 0.15 & 0.174 & 7.5 &  & 30 (73.2\%) & 83.10\\
  $50\%$ & 11.8 &  & 41 (100.0\%) & 0.21 & 0.288 & 8.5 &  & 41 (100.0\%) & 0.18 & 0.210 & 7.7 &  & 13 (31.7\%) & 200.64\\
  $n-1$ & 0.0 &  & 41 (100.0\%) & 0.28 & 0.362 & 8.4 &  & 41 (100.0\%) & 0.28 & 0.311 & 8.4 &  & 6 (14.6\%) & 678.36\\
  \midrule
  \multicolumn{15}{@{}l}{$70\le n_0\le100$, $20$ instances, $\card{\mathcal R}=10$}\\
  $1$ & 45.3 &  & 20 (100.0\%) & 1.24 & 2.243 & 12.7 &  & 20 (100.0\%) & 0.97 & 1.146 & 7.8 &  & 20 (100.0\%) & 3.05\\
  $10\%$ & 32.3 &  & 20 (100.0\%) & 1.84 & 2.917 & 7.6 &  & 20 (100.0\%) & 1.57 & 1.780 & 6.1 &  & 0 (0.0\%) & --\\
  $25\%$ & 21.5 &  & 20 (100.0\%) & 2.47 & 3.580 & 7.1 &  & 20 (100.0\%) & 2.22 & 2.477 & 6.3 &  & 0 (0.0\%) & --\\
  $50\%$ & 11.1 &  & 20 (100.0\%) & 3.14 & 4.428 & 7.1 &  & 20 (100.0\%) & 2.82 & 3.063 & 6.6 &  & 0 (0.0\%) & --\\
  $n-1$ & 0.0 &  & 20 (100.0\%) & 3.62 & 4.952 & 7.1 &  & 20 (100.0\%) & 3.54 & 3.809 & 7.1 &  & 0 (0.0\%) & --\\
  \bottomrule
\end{tabular}

}
\end{widetable}

{
These pair sets are easy for the compact formulation: with both
  penalty vectors, it solves every instance, with a geometric
mean total method time of at most $5.0$~seconds. The destination-aware
penalties strengthen the formulation when few destinations are
requested. With one destination per root, they are strictly smaller
than $M_{\texttt{rflow}}$ on an average of $38.2\%$ and $45.3\%$ of the
entries in the two blocks; the mean LP gap decreases from $13.8\%$ to
$9.1\%$ and from $12.7\%$ to $7.8\%$, respectively. Both compact
  configurations solve all instances with comparable
runtime performance; the destination-aware hybrid's main benefit is the
stronger relaxation when few destinations are requested. The
reduction in the LP gap shrinks as the number of destinations grows,
and it vanishes with all destinations, where the two penalty vectors
coincide by Proposition~\ref{prop:destination-hybrid}. The rule
$M_{\texttt{dtop}}$ is not tested, since $U_{\texttt{flow}}\le
U_{\texttt{top}}$ implies $M_{\texttt{dflow}}\le M_{\texttt{dtop}}$
componentwise.

The covering decomposition behaves differently. With one destination
per root, i.e., $5$ or $10$ requested pairs, it solves every instance.
Its geometric mean total times are $0.12$ and $3.05$~seconds in
the two blocks, compared with $0.122$ and $1.146$~seconds for the
hybrid, whose penalties, like the covering decomposition, only
involve the active roots. Thus, covering
is competitive on the smaller block, while the hybrid is faster on the
larger block. As the number of pairs grows, covering solves $41$, $30$,
$13$, and $6$ of the $41$ instances with $n_0\le60$ when
$10$--$30$, $20$--$75$, $40$--$150$, and $75$--$295$ pairs are
requested, respectively; it solves none of the $20$ larger instances
as soon as $70$ or more pairs are requested. The compact formulation
solves all of them. Larger pair sets can lead to larger unions of
selected paths and therefore weaker basic covering inequalities.
This provides a possible explanation for the observed deterioration in
covering performance.}

\subsection{The Arc-Destruction Setting of Israeli and Wood}
\label{sec:results-iw}
In the single-pair setting of~\citet{Israeli_Wood:2002}, we compare
their covering decomposition with the compact formulation on grid
networks generated as in their Section~7.1: an $r\times r$ grid with a
source joined to the first column and a sink joined to the last one,
arc lengths drawn uniformly in $\set{1,\ldots,10}$, unit interdiction
costs on the interdictable grid arcs, and budgets~$k\in\set{5,10}$, i.e., the $k$ most vital arcs
problem of their Table~5. We generate five instances for each grid and
budget. The source-to-first-column and last-column-to-sink arcs
are protected by assigning them interdiction cost~$k+1$, exceeding the
budget. The original networks are not available. Since our
implementation requires positive lengths and a strongly connected
graph, the arcs from the source and to the sink have length~1
instead of~0, which adds~$2$ to every source--sink distance, and
protected return arcs with length~$1$ and interdiction
cost~$k+1$ lead from every other vertex back to the source,
which no simple source--sink path uses. Neither change modifies the
optimal interdiction. These instances satisfy
Assumption~\ref{ass:feasibility}, since every source--sink cut has
cost larger than~$k$, but not Assumption~\ref{ass:strong-residual}:
removing the two arcs that enter a corner vertex of the grid makes it
unreachable from the source. The penalties are valid by
Propositions~\ref{prop:polynomial-rules}
and~\ref{prop:destination-hybrid}, which only require
Assumption~\ref{ass:feasibility}. The single pair makes
Formulation~\eqref{eq:single-level} a model with one dual block, and we
use the destination-aware penalties~$\min\set{M_{\texttt{rflow}},M_{\texttt{dflow}}}$, as in
\Cref{sec:results-mpsp}.

{
\Cref{tab:iw} reports the instances solved within $1{,}800$~seconds
and the geometric mean of the total time, including the
computation of the penalties.}

\begin{widetable}[htb]
  \centering
  \caption{Single-pair grids of~\citet[Table~5]{Israeli_Wood:2002}:
    instances solved to proven optimality (out of~$5$) and geometric
    mean of the total time (seconds) for the covering
    decomposition and for Formulation~\eqref{eq:single-level} with
    the destination-aware hybrid penalties~$\min\set{M_{\texttt{rflow}},M_{\texttt{dflow}}}$ (hybrid), without SVIs, with the
    cardinality (card.) SVIs, and with the cardinality SVIs and the
    objective cutoff. Time limit:
    $1{,}800$~seconds.}
  \label{tab:iw}
  \small
  \setlength{\tabcolsep}{3pt}
  \renewcommand{\arraystretch}{0.85}
  {
\begin{tabular}{@{}lr rr r rr r rr r rr@{}}
  \toprule
  & & \multicolumn{2}{c}{covering (basic SVIs)} & & \multicolumn{2}{c}{hybrid} & & \multicolumn{2}{c}{hybrid + card.\ SVIs} & & \multicolumn{2}{c}{hybrid + card.\ SVIs + cutoff}\\
  \cmidrule{3-4}\cmidrule{6-7}\cmidrule{9-10}\cmidrule{12-13}
  grid & $k$ & solved & time & & solved & time & & solved & time & & solved & time\\
  \midrule
  $7\times7$ & 5 & 5 & 0.2 &  & 5 & 0.5 &  & 5 & 0.4 &  & 5 & 0.5\\
  $7\times7$ & 10 & 5 & 0.3 &  & 5 & 8.0 &  & 5 & 8.6 &  & 5 & 16.8\\
  $8\times8$ & 5 & 5 & 0.4 &  & 5 & 0.8 &  & 5 & 0.8 &  & 5 & 0.7\\
  $8\times8$ & 10 & 5 & 0.8 &  & 5 & 25.2 &  & 5 & 10.1 &  & 5 & 36.2\\
  $9\times9$ & 5 & 5 & 0.9 &  & 5 & 1.4 &  & 5 & 1.0 &  & 5 & 1.3\\
  $9\times9$ & 10 & 5 & 1.4 &  & 5 & 40.7 &  & 5 & 17.5 &  & 5 & 86.8\\
  $12\times12$ & 5 & 5 & 4.2 &  & 5 & 6.8 &  & 5 & 4.9 &  & 5 & 5.5\\
  $12\times12$ & 10 & 5 & 6.8 &  & 5 & 373.1 &  & 5 & 62.2 &  & 5 & 979.3\\
  \bottomrule
\end{tabular}

}
\end{widetable}

{
The covering decomposition solves every instance and is the fastest
method in every grid--budget group, which confirms the ranking
of~\citet{Israeli_Wood:2002} with a modern solver. With the $1{,}800$-second time
limit, all compact configurations also solve all $40$ grids. The
objective cutoff still does not help in this setting: on the
$12\times12$ grids with $k=10$, the compact formulation with
cardinality SVIs has geometric mean time $62.2$ seconds,
whereas adding the cutoff raises it to $979.3$ seconds. In a separate
reproduction with the uniform delays of~\citet{Israeli_Wood:2002}, the
delay $d^+=10$ returned a suboptimal interdiction on one instance, so
that a uniform delay equal to the largest arc length is not valid in
general.}

The advantage of the covering decomposition has a structural
explanation. With a single pair, the optimal value is at least~$z$ if
and only if the interdiction hits every $s$--$t$ path shorter
than~$z$. This is a set-covering condition without lengths or
penalties. On the original grids with zero-length source and sink
connector arcs, the arithmetic mean gaps over the five instances for
each grid size and budget range from $0$ to $3.86\%$ for the covering-LP
bound. For the compact formulation with $k=10$, the corresponding means
range from $22.35\%$ to $33.00\%$, even with penalties computed from the
optimal value. \Cref{app:lp-diagnostic} defines both diagnostics and
reports the means for each configuration. With several pairs, the objective is a sum of
residual distances, and a lower bound on this sum is not a single
covering condition: an interdiction can increase some distances and
leave the others unchanged. The covering decomposition remains exact
for several pairs (\Cref{sec:svi}), but the support of each of its
inequalities is the union of the selected paths of all requested pairs,
so that each inequality excludes fewer interdictions as the number of
pairs grows, while the weakness of a uniform penalty grows with the
number of pairs as well.
{
The experiments of the previous sections confirm this behavior. With
$5$ or $10$ requested pairs, the covering decomposition and the compact
formulation both solve every instance of \Cref{tab:e3}. From a few tens
of pairs on, covering solves fewer and fewer instances, while the
compact formulation solves all of them. With all pairs, it solves $13$
of the $124$ scaling instances with $n_0\le50$ in \Cref{tab:e5}, against
$116$ for the compact formulation with $M_{\texttt{rflow}}$. Hence, on
the tested instances, the recommendation of~\citet{Israeli_Wood:2002}
to avoid penalties holds for a single pair on the tested grids, and the
compact formulation with strengthened penalties is the better choice
when many pairs are requested.}

\section{Conclusions and Future Work}
\label{sec:conclusions}

We studied the interdiction problem of the multiple-pairs shortest-path
problem, in which a single budget-constrained arc-removal decision
simultaneously affects a prescribed set~$\mathcal K$ of
origin-destination pairs. The all-pairs problem is recovered by
taking all ordered pairs of distinct vertices.
Our root-aggregated follower formulation remains totally unimodular:
for each active root only the vertices in~$\mathcal T_w$ have unit
demand, while the remaining non-root vertices are transit nodes. This gives an
exact compact MILP, whose strength depends on the penalties.

We characterized the valid penalties exactly through the residual dual
optimal face. With a proper destination subset, non-destination labels
can move while the requested destination labels stay fixed at their
residual distances, so that the penalties of different arcs are coupled
within a single interdiction. Over all feasible interdictions, however,
validity reduces to one threshold for each arc and root, for every pair
set: the componentwise smallest penalties~$\eta_{\mathcal K}(a,w)$
exist and only involve paths through a single interdicted arc toward
the requested destinations. In the all-pairs case they coincide
with~$\eta(a,w)$. We also derived the poly-time computable
destination-aware rules $M_{\texttt{dtop}}$ and $M_{\texttt{dflow}}$,
based on robust upper bounds for the requested targets.

Overall, the computational results favor the compact formulation with
the proposed surviving-path penalties, $M_{\texttt{rflow}}$ for
all pairs and its combination with $M_{\texttt{dflow}}$ for pair sets,
on the tested instances with up to $100$ vertices. It performs best
without supervalid inequalities in the all-pairs
experiments, among the tested lazy-SVI configurations.
Destination-aware hybrid penalties further strengthen the relaxation
when few destinations are requested. Among the tested settings, the covering decomposition
of~\citet{Israeli_Wood:2002} is preferable only for single-pair
interdiction on grids and remains competitive when a few pairs are
requested.

Further extensions include the max-min version over the requested
pairs, blocker counterparts, weighted pair objectives, and pair sets
that change over time.

\section*{Acknowledgements}

The authors thank Valerio Dose and Yanlu Zhao for their contributions
to earlier stages of this work and for many fruitful discussions.

\section*{Declaration of AI Use}

The authors used generative AI tools, primarily Claude (Anthropic),
and additionally ChatGPT and Codex (OpenAI), between September and
November 2026 to assist with software development, code verification,
computational data analysis, and the drafting, editing, and checking of
parts of the manuscript, including mathematical arguments and proofs.

AI-assisted outputs were critically reviewed by the authors, and the
computational implementations were tested against independently computed
reference solutions. The authors retained full control over the
scientific contributions and take full responsibility for the accuracy,
correctness, and integrity of the work.

\bmsubsection*{Data Availability Statement}
The instances, the source code of the solver, and the results of all
computational experiments are available at \url{https://github.com/AnnaLivia/MPSP-Interdiction-Problem} and \url{https://doi.org/10.5281/zenodo.24680030}.

\printbibliography

\clearpage
\appendix
\markboth{\MakeUppercase{\shortauthors}}{\MakeUppercase{\shorttitle}}

\section{Auxiliary Single-Pair LP Diagnostics}
\label{app:lp-diagnostic}

The diagnostics in \Cref{sec:results-iw} use the original grids
with zero-length source and sink connector arcs. Let $z^*$ be the
optimal interdicted $s$--$t$ distance under budget~$k$. For
$L\in\Z_{\ge0}$, let $C(L)$ be the optimal value of the linear
optimization problem that minimizes $\sum_{a\in\mathcal A}x_a$ subject
to $0\le x_a\le1$ for every arc~$a$, $x_a=0$ on protected arcs, and
$\sum_{a\in P}x_a\ge1$ for every $s$--$t$ path~$P$ with
$\ell(P)\le L$. This is a minimum fractional path-cover problem,
distinct from the cost-maximizing master in \Cref{sec:svi}. The
reported covering-LP distance upper bound is
\[
  z^{\mathrm{cov}}\define
  \min\Set{L\in\Z_{\ge0}:C(L)>k}.
\]
Indeed, $C(L)>k$ implies that no feasible interdiction can hit
every $s$--$t$ path of length at most~$L$, so $z^*\le L$.

The compact diagnostic uses the destination rule with
$U(s,t)=z^*$. For each interdictable arc~$a=(u,v)$, it sets
\[
  M_{\mathrm{oracle}}(a,s)\define
  \max\Set{0,\ z^*-d(s,u)-\ell_a-d(v,t)}.
\]
As $z^*$ bounds every feasible interdicted $s$--$t$ distance,
the validity argument of \Cref{prop:destination-polynomial} applies.
These oracle penalties are not asserted to equal the exact
thresholds~$\eta_{\mathcal K}(a,s)$.

The gaps in \Cref{tab:lp-diagnostic} are $100(B-z^*)/z^*$,
where~$B$ is either $z^{\mathrm{cov}}$ or the compact linear
optimization relaxation value under the oracle penalties.

\begin{table}[htb]
  \centering
  \caption{Arithmetic mean relative gaps (\%) on the original
    grids with zero-length connector arcs, over five instances per
    grid size and budget.}
  \label{tab:lp-diagnostic}

\begin{tabular}{@{}crrr@{}}
  \toprule
  grid & $k$ & covering-LP bound & oracle compact LP\\
  \midrule
  $7\times7$   &  5 & 1.25 & 19.15\\
  $7\times7$   & 10 & 3.86 & 33.00\\
  $8\times8$   &  5 & 0.00 & 14.64\\
  $8\times8$   & 10 & 1.60 & 30.16\\
  $9\times9$   &  5 & 0.00 & 11.64\\
  $9\times9$   & 10 & 1.46 & 27.24\\
  $12\times12$ &  5 & 0.69 &  9.55\\
  $12\times12$ & 10 & 3.14 & 22.35\\
  \bottomrule
\end{tabular}

\end{table}

\end{document}